\documentclass[11pt,a4paper]{article}

\usepackage[margin=1in]{geometry}
\usepackage{amsmath,amssymb,amsthm,mathtools}
\usepackage{enumitem}
\usepackage{hyperref}
\usepackage{mathrsfs}
\usepackage{authblk}
\usepackage[nocompress]{cite}
\usepackage{bm}
\usepackage{tikz-cd}

\hypersetup{
  colorlinks=true,
  linkcolor=blue,
  citecolor=blue,
  urlcolor=blue
}

\newtheorem{theorem}{Theorem}[section]
\newtheorem{proposition}[theorem]{Proposition}
\newtheorem{lemma}[theorem]{Lemma}
\newtheorem{corollary}[theorem]{Corollary}
\newtheorem{example}[theorem]{Example}
\theoremstyle{definition}
\newtheorem{definition}[theorem]{Definition}
\newtheorem{remark}[theorem]{Remark}

\title{Invariant complex Finsler metrics on reductive complex homogeneous spaces}
\author{%
  Kuankuan Luo$^{1,\ast}$ \qquad Chunping Zhong$^{1,\dagger}$\\
  {\small $^{1}$School of Mathematical Sciences, Xiamen University,
  Xiamen 361005, P.~R.~China}\\
  {\small $^{\ast}$E-mail: \texttt{luokk2026@163.com}}\\
  {\small $^{\dagger}$Corresponding author. E-mail: \texttt{zcp@xmu.edu.cn}}
}
\date{}

\usepackage{xcolor}

\begin{document}

\maketitle

\begin{abstract}
 Let $(M=G/H,J)$ be a reductive complex homogeneous space and $F$ be a $G$-invariant strongly pseudoconvex complex Finsler metric on $M$. Using the adapted $(1,0)$-frame field obtained by parallel transformation with respect to the canonical connection on $M$, we obtain the necessary and sufficient condition for $F$ to be a K\"ahler-Finsler metric (resp. weakly K\"ahler-Finsler metric) in terms of Lie-bracket. The curvature tensor of the complex Rund connection associated to $F$ is then decomposed into an algebraic part and a derivative part, yielding explicit formulas for the holomorphic bisectional curvature and holomorphic sectional curvature of $F$, respectively. As an application,  we give a Lie-theoretic proof of the rigidity theorem: every $G$-invariant strongly pseudoconvex complex Finsler metric on a complex symmetric manifold $M=G/H$ is necessarily a K\"ahler--Berwald metric.
\end{abstract}


\tableofcontents

\section{Introduction and main results}\label{sec:1}
\hspace{0.4cm}
%A basic principle of differential geometry on homogeneous Riemannian manifolds is that invariant geometric structures can be studied through Lie group and Lie algebraic method at a single point. 
On a reductive homogeneous space $M=G/H$, a $G$-invariant Riemannian metric corresponds to an $\operatorname{Ad}(H)$-invariant inner product on $T_oM\cong\mathfrak m$. This correspondence underlies the classical theory of homogeneous Riemannian manifolds \cite{Wallach,Heintze}. From the viewpoint of  Chern \cite{ChernSSq}, Finsler geometry is Riemannian geometry without the quadratic restriction. In  \cite{DengHou2004},  Deng--Hou developed an algebraic description of invariant real Finsler metrics using Minkowski Lie pairs. Huang \cite{Huang2015,Huang2017} subsequently introduced the key spray vector and connection operator to express the flag, Ricci, Landsberg and $S$-curvatures of homogeneous real Finsler spaces in terms of Lie algebra data. For more related results on the curvatures of homogeneous real Finsler spaces, we refer to \cite{DengHu2013,YanDeng2014,XuDeng2017,HuangShen2019} and references therein.

The algebraic study of homogeneous complex Finsler spaces was first carried out in Deng--Hou \cite{DengHou2006}, see also Deng \cite{Deng2012}. They described homogeneous complex Finsler spaces as coset spaces and used Minkowski representations of Lie groups and Lie algebras to study invariant complex Finsler metrics on complex homogeneous spaces.

 It is well-known that the Bergman metric \cite{Bergmann1922} plays an important role in the geometric function theory of several complex variables and complex geometry. In complex Finsler geometry, experts concern whether there exists holomorphically invariant non-Hermitian quadratic strongly pseudoconvex Finsler metrics (cf. \cite{BK}). In \cite{Zhong2023}, Zhong explicitly constructed a class of non-Hermitian quadratic holomorphically invariant strongly convex  K\"ahler-Berwald metrics on the unit polydisk in $\mathbb{C}^n(n\geq 2)$,  and then in \cite{GeZhong2024}, Ge and Zhong explicitly constructed a class of non-Hermitian quadratic  holomorphically invariant K\"ahler--Berwald metrics on the irreducible bounded symmetric domains of type $I-IV$ with rank $\geq 2$. The complete characterizations were immediately obtained on the unit polydisk \cite{LinZhong2023},  polyball \cite{LinWangZhong2024}, the irreducible bounded symmetric domains of type $I-IV$ \cite{Zhong2025}, and the complex Grassmannian manifold \cite{CaoGeZhong2024}. The corresponding Schwarz lemma for holomorphic mappings with respect to these invariant metrics were also established on the polydisk \cite{LinZhong2023} and the irreducible bounded symmetric domains of type $I-IV$ \cite{Zhong2025}, which generalized the classical Schwarz lemma previously established by Lu \cite{Look1957,Look1958}  with respect to the Bergman metric.

Let $(M=G/H,J)$ be a reductive complex homogeneous space, where $G$ is a connected real Lie group acting transitively by biholomorphisms, and fix a reductive decomposition
$$
\mathfrak g=\mathfrak h\oplus\mathfrak m,\qquad \operatorname{Ad}(H)\mathfrak m\subset\mathfrak m.
$$
Recall that $G/H$ is reductive in each of the following cases \cite{Nomizu_K}:
\begin{enumerate}
\item[(1)] $H$ is compact;
\item[(2)] $H$ is connected and $\mathfrak h$ is reductive in $\mathfrak g$, meaning that the adjoint representation of $\mathfrak h$ on $\mathfrak g$ is completely reducible. In particular, this holds when $H$ is connected and semisimple;
\item[(3)] $H$ is discrete, in which case $\mathfrak h=0$ and $\mathfrak m=\mathfrak g$.
\end{enumerate}

In this paper, we study $G$-invariant strongly pseudoconvex complex Finsler metrics on $(G/H,J)$, extending our earlier results \cite{LuoXiaoZhong2025}  for left-invariant complex Finsler metrics on Lie groups (which are simultaneously complex manifolds) to reductive complex homogeneous spaces.

Compared with the Lie group case \cite{LuoXiaoZhong2025}, the main difficulty is that the structure functions of the adapted $(1,0)$-frame on $G/H$ generally vary with the base point $z\in \widetilde{\mathcal{U}}$, where $\widetilde{\mathcal{U}}$ is a neighborhood of $o$. To compute the curvature at $(o;v)\in T^{1,0}M$ with $v\neq 0$, we need the values of the adapted-frame structure functions $\lambda_{\;j\bar{k}}^i(o), \lambda_{\;j\bar{k}}^{\bar{i}}(o)$ and  first directional derivatives $e_l(\lambda_{\;j\bar{k}}^i)(o),e_{\bar{l}}(\lambda_{\;j\bar{k}}^i)(o)$. These derivatives are taken along the adapted $(1,0)$- and $(0,1)$-frame fields on the base manifold, respectively. We determine these data by applying Lemma~\ref{lm:3.9} to the canonical connection $\overset{\mathbf{C}}{\nabla}$  on $G/H$.

To present our main results, let $g_v$ be the fundamental tensor,  $\mathcal C_v^+$ and $\mathcal{C}_v^-$ be  the complex Cartan tensors associated to $F$ respectively, $\mathcal S_v$ be the symmetric product  and $\overline{\mathcal S}_v$ its complex conjugate introduced in Section~\ref{sec:3.1}. We also let $\mathcal N_v$  be the connection operator (cf. Definition~\ref{df:3.14}) and $\sigma:\mathfrak{g}^{\mathbb{C}}\to\mathfrak{g}^{\mathbb{C}}$  be the conjugation with respect to the real form $\mathfrak{g}$, namely
\begin{equation*}
    \sigma(X+\sqrt{-1}Y)=X-\sqrt{-1}Y,\qquad X,Y\in \mathfrak{g}.
\end{equation*}

\begin{theorem}[cf. Theorem~\ref{th:3.16}]\label{th:1.1}
Let $(M=G/H,J)$ be a reductive complex homogeneous space with reductive decomposition $\mathfrak g=\mathfrak h\oplus\mathfrak m$ and $\operatorname{Ad}(H)\mathfrak m\subset\mathfrak m$. Let $F:T^{1,0}M\to[0,\infty)$ be a $G$-invariant strongly pseudoconvex complex Finsler metric. Then
\begin{enumerate}
\item[(1)] $F$ is a K\"ahler--Finsler metric if and only if for every $v\in\mathfrak m^{1,0}\setminus\{0\}$ and all $u,w\in\mathfrak m^{1,0}$,
\begin{align*}
&g_v(w,[u,\sigma(v)]_{\mathfrak m^{1,0}})-g_v([v,w]_{\mathfrak m^{1,0}},u)-g_v(v,[u,\sigma(w)]_{\mathfrak m^{1,0}})\nonumber\\
&-\overline{\mathcal{S}}_v(u,[v,\sigma(w)]_{\mathfrak m^{1,0}})-\mathcal C_v^+(w,u,\mathcal N_v(v))+\mathcal C_v^-(w,u,[v,\sigma(v)]_{\mathfrak m^{1,0}})=0;
\end{align*}
\item[(2)] $F$ is a weakly K\"ahler--Finsler metric if and only if for every $v\in\mathfrak m^{1,0}\setminus\{0\}$ and all $w\in\mathfrak m^{1,0}$,
\begin{align*}
g_{v}(w,[v,\sigma(v)]_{\mathfrak{m}^{1,0}})-g_{v}([v,w]_{\mathfrak{m}^{1,0}},v)-g_{v}(v,[v,\sigma(w)]_{\mathfrak{m}^{1,0}})-\mathcal{S}_{v}(\mathcal N_v(v),w)=0.
\end{align*}
\end{enumerate}
\end{theorem}

\begin{remark}
If in particular $H=\{e\}$, the space $(G/H,J)$ reduces to the Lie group $G$ (which is simultaneously a complex manifold) equipped with a $G$-invariant complex structure, and $G$-invariant complex Finsler metrics are precisely left-invariant ones. In this case, the above theorem recovers the corresponding characterizations of K\"ahler--Finsler metrics (resp. weakly K\"ahler--Finsler metrics) obtained in \cite{LuoXiaoZhong2025}.
\end{remark}

Let $\mathbf R_v$ denote the horizontal $(1,1)$-curvature operator of the complex Rund connection associated to $F$. In the decomposition $\mathbf R_v=\mathfrak R_v+\mathcal R_v$, $\mathfrak R_v$ represents the contribution involving derivatives of the adapted-frame structure functions, while $\mathcal R_v$ represents the remaining algebraic contribution; see \eqref{eq:3.13} and \eqref{eq:3.14}, respectively. We use the standard flat connection $\mathcal D$ on $\mathfrak m^{1,0}\setminus\{0\}$ to differentiate fiber-dependent quantities in both holomorphic and antiholomorphic directions, as specified in \eqref{eq:3.40}.

Our second main result gives explicit formulas for the holomorphic bisectional curvature  of $(G/H,J,F)$.

\begin{theorem}[cf. Theorem~\ref{th:3.18}]\label{th:1.3}
Let $(M=G/H,J)$ be a reductive complex homogeneous space with reductive decomposition $\mathfrak g=\mathfrak h\oplus\mathfrak m$ and $\operatorname{Ad}(H)\mathfrak m\subset\mathfrak m$. Let $F:T^{1,0}M\to[0,\infty)$ be a $G$-invariant strongly pseudoconvex complex Finsler metric. Let $v\in\mathfrak{m}^{1,0}\setminus\{0\}$ be the radial vector with variable direction and $w\in\mathfrak{m}^{1,0}$ be an arbitrary vector. Then for the curvature tensor decomposition $\mathbf{R}^i_{\;k\bar{j}}=\mathfrak{R}^i_{\;k\bar{j}}+\mathcal{R}^i_{\;k\bar{j}}$ as in \eqref{eq:3.12}, we have
\begin{equation*}
    \begin{split}
g_v(\mathbf{R}_v(w,\sigma(w))v, v)=g_v([[\sigma(w), w]_{\mathfrak{h}^{\mathbb{C}}}, v]_{\mathfrak{m}^{1,0}}, v)+g_v\left(\mathcal{R}_v(w,\sigma(w))v,v\right),
    \end{split}
\end{equation*}
where $\mathcal{R}_v(w,\sigma(w))v$ is given by
\begin{equation*}
	\begin{split}
	\mathcal{R}_v(w, \sigma(w))v=&\mathcal{D}_{\sigma(\mathcal N_v(w))}\mathcal N_v(w)-\mathcal{D}_{[v, \sigma(w)]_{\mathfrak{m}^{1, 0}}}\mathcal N_v(w)-\mathcal N_v([w, \sigma(w)]_{\mathfrak{m}^{1, 0}})\\
    &+[\mathcal N_v(w), \sigma (w)]_{\mathfrak{m}^{1, 0}}-[v, [\sigma(w), w]_{\mathfrak{m}^{0,1}}]_{\mathfrak{m}^{1, 0}}.
	\end{split}
\end{equation*}
\end{theorem}

\begin{remark}
If $H=\{e\}$, then $\mathfrak R_v=0$ and
$$
\mathbf R_v(w,\sigma(w))v=\mathcal R_v(w,\sigma(w))v.
$$
Thus formula  for $g_v(\mathbf{R}_v(w,\sigma(w))v, v)$ given in Theorem \ref{th:1.3} recovers the holomorphic bisectional curvature obtained in \cite{LuoXiaoZhong2025}.
\end{remark}

Finally, if $(G/H, J)$ is a complex symmetric coset space, then its canonical decomposition satisfies $[\mathfrak m,\mathfrak m]\subset\mathfrak h$. By Theorem \ref{th:1.1} and the formula of the complex Rund connection coefficients, we are able to obtain a rigidity theorem: every $G$-invariant strongly pseudoconvex complex Finsler metric $F$ on such a space is  necessarily a K\"ahler--Berwald metric (Theorem~\ref{th:4.5}).
\begin{theorem}\label{th:1.6}
Let $(M,J,F)$ be a globally symmetric complex Finsler space. Then $F$ is necessarily a K\"ahler--Berwald metric.
\end{theorem}
\begin{remark}
Theorem \ref{th:1.6} was also obtained by Ge--Zhou \cite{GeZhou} and He--Zhong \cite{HeZhong} independently. Our proof in this paper is based Lie theory which is essentially different from both \cite{GeZhou} and \cite{HeZhong}.
\end{remark}

In the symmetric case, the holomorphic bisectional curvature and holomorphic sectional curvature formulas simplify to the following Lie-algebraic expressions.

\begin{theorem}[cf. Theorem~\ref{th:4.7}]
Let $(M,J,F)$ be a globally symmetric complex Finsler space. Identify $M$ with a complex symmetric coset space $(G/H,J,\Theta)$ as in Theorem~\ref{th:4.2}, and let $\mathfrak g=\mathfrak h\oplus\mathfrak m$ be the canonical decomposition. Let $v,w\in\mathfrak m^{1,0}\setminus\{0\}$ satisfying $g_v(v,v)=g_v(w,w)=1$. Then the holomorphic bisectional curvature in the directions $v$ and $w$, and the holomorphic sectional curvature in the direction $v$, are given respectively by
\begin{align*}
B(v,w)=g_v\bigl([[\sigma(w),w]_{\mathfrak h^{\mathbb C}},v]_{\mathfrak m^{1,0}},v\bigr),\quad K(v)=2g_v\bigl([[\sigma(v),v]_{\mathfrak h^{\mathbb C}},v]_{\mathfrak m^{1,0}},v\bigr).
\end{align*}
\end{theorem}

This paper is organized as follows. In Section~\ref{sec:2}, we recall the background on reductive complex homogeneous spaces and the complex Rund connection. In Section~\ref{sec:3}, we construct the adapted frame and derive the K\"ahler--Finlser (resp. weakly K\"ahler--Finsler) criteria and curvature formulas. In Section~\ref{sec:4}, we apply these formulas to give a characterization of globally symmetric complex Finsler spaces.

\section{Preliminaries}\label{sec:2}
\hspace{0.4cm}
In this section we recall some basic facts and fix some notation about reductive complex homogeneous spaces and complex Finsler geometry. For more details, we refer to \cite{LuoXiaoZhong2025,FoundationOfDG2,bao2012introduction,AbateAndPatrizio}.

\subsection{Reductive complex homogeneous spaces}\label{sec:2.1}
\hspace{0.4cm}
Let $M=G/H$ be a homogeneous space, where $G$ is a connected real Lie group and $H$ is the closed isotropy subgroup at $o=eH$, $\pi:G\rightarrow G/H$ the natural projection.  Let $\mathfrak g=\operatorname{Lie}G$ and $\mathfrak h=\operatorname{Lie}H$. The homogeneous space $G/H$ is called reductive if there exists a vector subspace $\mathfrak m\subset\mathfrak g$ such that
$$
\mathfrak g=\mathfrak h\oplus\mathfrak m, \qquad \operatorname{Ad}(H)\mathfrak m\subset\mathfrak m.
$$
In what follows, we assume that $G/H$ is reductive and fix such a decomposition. The restriction
$$
d\pi_e|_{\mathfrak m}:\mathfrak m\longrightarrow T_oM
$$
is then a vector space isomorphism, which we denote simply by $d\pi$.

Suppose that $M$ is equipped with a $G$-invariant almost complex structure $J$. Under the identification above, its value at $o$ determines a linear operator
$$
I:=(d\pi)^{-1}\circ J_o\circ d\pi
$$
on $\mathfrak m$. It satisfies
$$
I^2=-\operatorname{Id}_{\mathfrak m}, \qquad I\circ\operatorname{Ad}(h)=\operatorname{Ad}(h)\circ I, \qquad h\in H.
$$
If $H$ is connected, the second of the above conditions is equivalent to
$$
I\circ \operatorname{ad}(X)=\operatorname{ad}(X)\circ I,\qquad \forall X\in\mathfrak h.
$$
Assume $J$ is integrable, then the integrability condition for $J$ implies
\begin{equation*}
[IX,IY]_{\mathfrak m}-[X,Y]_{\mathfrak m}-I[X,IY]_{\mathfrak m}-I[IX,Y]_{\mathfrak m}=0,\qquad \forall X,Y\in\mathfrak m.
\end{equation*}
Here and below $[\cdot,\cdot]_{\mathfrak m}$ denotes the projection onto $\mathfrak m$. In this case, we call $(M,J)$ a reductive complex homogeneous space.
Throughout the rest of this paper, unless otherwise stated, every reductive homogeneous space is assumed to be complex in this sense.

Let $\mathfrak m^{\mathbb C}=\mathfrak m\otimes_{\mathbb R}\mathbb C$. The natural complex-linear extension of $I$ to $\mathfrak{m}^{\mathbb{C}}$ induces the following decomposition
$$
\mathfrak m^{\mathbb C}=\mathfrak m^{1,0}\oplus\mathfrak m^{0,1},
$$
where $\mathfrak m^{1,0}$ and $\mathfrak m^{0,1}$ are the $+\sqrt{-1}$ and $-\sqrt{-1}$ eigenspaces of $I$, respectively. Equivalently,
$$
\mathfrak m^{1,0}=\{X-\sqrt{-1}IX:X\in\mathfrak m\},\qquad \mathfrak m^{0,1}=\{X+\sqrt{-1}IX:X\in\mathfrak m\}.
$$
The complex-linear extension of $d\pi$ gives a complex isomorphism
$$
d\pi:\mathfrak m^{1,0}\longrightarrow T_o^{1,0}M.
$$

Let $\{E_1,\dots,E_m,IE_1,\dots,IE_m\}$ be a real basis of $\mathfrak m$, and set
$$
\mathcal E_i:=E_i-\sqrt{-1}\,IE_i\in\mathfrak m^{1,0}.
$$
Then $\{\mathcal E_i\}$ is a basis of $\mathfrak m^{1,0}$, which satisfies $I(\mathcal{E}_i)=\sqrt{-1}\mathcal{E}_i$. Let $\{\mathcal F_\alpha\}$ be a basis of $\mathfrak h$. Define
$\widetilde I:\mathfrak g\to\mathfrak g$ by setting
$$
\widetilde I\vert_{\mathfrak h}=0,\qquad \widetilde I\vert_{\mathfrak m}=I,
$$
and then extend it complex-linearly to $\mathfrak g^{\mathbb C}$. The same set $\{\mathcal F_\alpha\}$ serves as a complex basis of $\mathfrak h^{\mathbb C}$. Thus
\begin{equation}
\mathfrak g^{\mathbb C}=\mathfrak h^{\mathbb C}\oplus\mathfrak m^{\mathbb C}.\label{eq:2.1}
\end{equation}

Let $\sigma:\mathfrak g^{\mathbb C}\to\mathfrak g^{\mathbb C}$ be the conjugation with respect to the real form $\mathfrak g$, which is defined by
$$
\sigma(X+\sqrt{-1}Y)=X-\sqrt{-1}Y,\qquad X,Y\in\mathfrak g.
$$
It is easy to check that $\sigma$ preserves each factors of \eqref{eq:2.1}. In particular,
$\sigma(\mathcal E_i)\in \mathfrak m^{0,1}$ and $\mathfrak m^{0,1}=\sigma(\mathfrak m^{1,0})$.

For $h\in H$, the induced action on $\mathfrak m^{1,0}$ is defined by
$$
\operatorname{Ad}_{\mathfrak m^{1,0}}(h)\mathcal E:=\operatorname{Ad}(h)E-\sqrt{-1}I(\operatorname{Ad}(h)E),\quad \mathcal E=E-\sqrt{-1}IE\in\mathfrak m^{1,0}.
$$
Its infinitesimal version is denoted by $\operatorname{ad}_{\mathfrak m^{1,0}}(\mathfrak h)$.

The integrability condition and the mixed Lie brackets of $\mathfrak{m}^{\mathbb{C}}$ imply
$$
[\mathfrak m^{1,0},\mathfrak m^{1,0}]\subset \mathfrak h^{\mathbb C}\oplus\mathfrak m^{1,0},\qquad [\mathfrak m^{1,0},\mathfrak m^{0,1}]\subset \mathfrak h^{\mathbb C}\oplus\mathfrak m^{1,0}\oplus\mathfrak m^{0,1}.
$$
We write the complexified Lie brackets in terms of structure constants as
$$
[\mathcal E_j,\mathcal E_k]=c^i_{\;jk}\,\mathcal E_i+c^\alpha_{\;jk}\,\mathcal F_\alpha,\quad
[\mathcal E_j,\sigma(\mathcal E_k)]=c^i_{\;j\bar k}\,\mathcal E_i+c^{\bar i}_{\;j\bar k}\,\sigma(\mathcal E_i)+c^\alpha_{\;j\bar k}\,\mathcal F_\alpha.
$$
Here the coefficients $c^i_{\;jk}$, $c^i_{\;j\bar k}$, $c^{\bar i}_{\;j\bar k}$, $c^\alpha_{\;jk}$, $c^\alpha_{\;j\bar k}$ are complex numbers. Indices $i,j,k,l,p,q,\dots$ run from $1$ to $m$, while $\alpha,\beta$ run from $1$ to $\dim_{\mathbb{R}}\mathfrak{h}$.

Applying $\sigma$ to the above bracket relations and comparing coefficients gives the conjugation relations
$$
\overline{c^i_{\;jk}}=c^{\bar i}_{\;\bar j\bar k},\qquad \overline{c^i_{\;j\bar k}}=c^{\bar i}_{\;\bar j k},\qquad \overline{c^\alpha_{\;j\bar k}}=c^\alpha_{\;\bar j k},
$$
where the over line on the left hand side of the above equalities denotes conjugation of complex numbers. Similarly, since $[\mathfrak{h},\mathfrak{m}]\subset \mathfrak{m}$, the $\mathfrak{h}$-component of $[\mathfrak{h},\mathfrak{m}]$ vanishes. It follows that
$$
[\mathcal F_\alpha,\mathcal E_j]=c^k_{\;\alpha j}\,\mathcal E_k\qquad \text{and}\qquad \overline{c^k_{\;\alpha j}}=c^{\bar k}_{\;\alpha\bar j}.
$$
Let $v=v^i\mathcal E_i\in\mathfrak m^{1,0}$ be a nonzero vector. We use the index $o$ (resp. $\bar o$) to denote contraction with $v^k$ (resp. $\bar v^k$), e.g.
$$
c^i_{\;jo}:=c^i_{\;jk}v^k,\qquad c^i_{\;j\bar{o}}:=c^i_{\;j\bar{k}}\bar{v}^k.
$$

\subsection{Complex Finsler metrics and curvatures}
\hspace{0.4cm}
A complex Minkowski norm on a complex vector space $V$ is a continuous function $F:V\to[0,\infty)$ which is smooth and positive on $V\setminus\{0\}$ and satisfies $F(\lambda v)=|\lambda|F(v)$ for all $\lambda\in\mathbb C$ and $v\in V$. The pair $(V,F)$ is called a complex Minkowski space. Given a basis $\{\xi_i\}$ and writing $v=v^i\xi_i$, the Levi matrix of $F$ is defined by
$$
g_{i\bar j}(v):=\frac{\partial^2F^2}{\partial v^i\partial\bar v^j}.
$$

A complex Finsler metric on a complex manifold $M$ is a continuous function $F:T^{1,0}M\to[0,\infty)$ which is smooth on $T^{1,0}M\setminus\{0\}$ and whose restriction to each fiber $T_z^{1,0}M$ is a complex Minkowski norm. The metric is called strongly pseudoconvex if the Levi matrix
$$
\left(g_{i\bar j}(z;v)\right):=\left(\frac{\partial^2F^2}{\partial v^i\partial\bar v^j}\right)
$$
is positive definite for all $z\in M$ and all nonzero vectors $v\in T_z^{1,0}M$. The coefficients $C_{i\bar j k}$ and $C_{i\bar j\bar k}$ are the complex Cartan tensors of $F$, defined by
$$
C_{i\bar j k}:=\frac{\partial g_{i\bar j}}{\partial v^k},\qquad C_{i\bar j\bar k}:=\frac{\partial g_{i\bar j}}{\partial\bar v^k}.
$$

Unless stated otherwise, all complex Finsler metrics in this paper are assumed to be strongly pseudoconvex. We shall use the fundamental tensor $g_{i\bar j}$ and its inverse $g^{\bar j i}$ to raise and lower indices; for example, given a structure constant $c^l_{\;pk}$, we write $c^{\;\,\bar i}_{\bar j\;k}=g_{l\bar j}c^l_{\;pk}g^{\bar i p}$, and similarly for other index conversions.

Let $\tilde M=T^{1,0}M\setminus\{0\}$ and $\tilde{\pi}:\tilde M\to M$ be the natural projection. Choose a smooth local frame $\{e_i\}$ of $T^{1,0}M$, with dual coframe $\{\varphi^i\}$, and use the same notation for their pullbacks. Write $v=v^ie_i$.

\begin{theorem}[{\cite[Theorem~3.1]{LuoXiaoZhong2025}}]
The bundle $\tilde{\pi}^\ast T^{1,0}M$ admits a unique complex Rund connection whose connection $1$-forms $\{\omega^i_j\}$ are characterized by
$$
d\varphi^i=\varphi^j\wedge\omega^i_j+\tau^i,
$$
where $\tau^i$ contain only horizontal $(2,0)$-components, and
$$
dg_{i\bar j}=g_{l\bar j}\,\omega^l_i+g_{i\bar l}\,\omega^{\bar l}_{\bar j}+C_{i\bar j k}\,\delta v^k+C_{i\bar j\bar k}\,\delta\bar v^k,\qquad\delta v^i=dv^i+v^j\omega^i_j.
$$
\end{theorem}
More precisely, we have
$$
\omega^i_j=\Gamma^i_{\;jk}\varphi^k+\Gamma^i_{\;j\bar k}\bar\varphi^k,\qquad\tau^i=T^i_{\;jk}\varphi^j\wedge\varphi^k
$$
and
$$
N^i_{\;k}:=v^j\Gamma^i_{\;jk},\qquad N^i_{\;\bar k}:=v^j\Gamma^i_{\;j\bar k},
$$
so that
$$
\delta v^i =dv^i+N^i_{\;k}\varphi^k +N^i_{\;\bar k}\bar\varphi^k.
$$
The forms $\{\varphi^i,\bar\varphi^i, \delta v^i,\delta\bar v^i\}$ constitute the associated horizontal and vertical coframe on $\tilde{M}$.

The complex Rund connection determines the horizontal subspace at $(z;v)\in\tilde{M}$, namely
$$
(\mathcal H_{\mathbb C})_{(z;v)}=\left\{X\in (T_{\mathbb C}\tilde M)_{(z;v)}:\delta v^i(X)=\delta\bar v^i(X)=0,\ \forall i\right\}.
$$
The restriction of $d\tilde\pi_{(z;v)}$ to this subspace is an isomorphism onto the fiber space $(T_{\mathbb C}M)_z$. We denote its inverse, the horizontal lift, by
$$
\chi_{(z;v)}:(T_{\mathbb C}M)_z\longrightarrow(\mathcal H_{\mathbb C})_{(z;v)}.
$$
In particular, horizontal lifting commutes with complex conjugation. These constructions are independent of the choice of pulled-back smooth frames; see \cite[Section~3.2]{LuoXiaoZhong2025}.

\begin{definition}[\cite{AbateAndPatrizio},\cite{Aikou-a}]
The complex Finsler metric $F$ is called K\"ahler--Finsler if $T^i_{\;jk}v^j=0$; called weakly K\"ahler--Finsler if $g_{i\bar l}\,T^i_{\;jk}v^j\bar v^l=0$; called  complex Berwald if $\Gamma^i_{\;jk}$ are independent of $v$ and $\bar v$;  called K\"ahler--Berwald if it is both K\"ahler--Finsler and complex Berwald.
\end{definition}

Let $\overset{\textbf{CR}}{\nabla}$ denote the complex Rund connection on $\tilde\pi^*T^{1,0}M$ associated to $F$. Its curvature is
$$
\Omega(X,Y)\xi=\overset{\textbf{CR}}{\nabla}_X\overset{\textbf{CR}}{\nabla}_Y\xi -\overset{\textbf{CR}}{\nabla}_Y\overset{\textbf{CR}}{\nabla}_X\xi -\overset{\textbf{CR}}{\nabla}_{[X,Y]}\xi,
$$
where $X,Y$ are complex vector fields on $\tilde M$ and $\xi$ is a smooth section of $\tilde\pi^*T^{1,0}M$. The curvature forms, defined by $\Omega(X,Y)e_j=\Omega^i_j(X,Y)e_i$, are
$$
\Omega^i_j=d\omega^i_j-\omega^k_j\wedge\omega^i_k.
$$

Let $R^i_{\;jk\bar l}$ denote the coefficient of $\varphi^k\wedge\bar\varphi^l$ in $\Omega^i_j$ with respect to the associated coframe. Equivalently, for $u,w\in T_z^{1,0}M$,
$$
\Omega_{(z;v)}\bigl(\chi_{(z;v)}(w),\chi_{(z;v)}(\bar w)\bigr)u=R^i_{\;jk\bar l}(z;v)\,u^jw^k\bar w^l e_i,
$$
where $u$ is identified with its value in the pullback fiber at $(z;v)$. Thus the coefficients used below evaluate the curvature on horizontal lifts, with the last argument in the pullback bundle.  For $v\ne0$, put $g_v(u,w)=g_{i\bar j}(z;v)\,u^i\bar w^j$.

\begin{definition}
For nonzero $v,w\in T_z^{1,0}M$, the holomorphic bisectional curvature with reference direction $v$ and the holomorphic sectional curvature are defined by
$$
B(v,w) =\frac{g_{i\bar j}R^i_{\;pk\bar l}\,v^p\bar v^j w^k\bar w^l}{g_v(v,v)\,g_v(w,w)},\qquad K(v)=2B(v,v).
$$
\end{definition}

We shall also use the contracted curvature forms
\begin{equation}\label{eq:2.2}
\begin{aligned}
\Omega^i&:=d(\delta v^i)-\delta v^j\wedge\omega^i_j =v^j\Omega^i_j\\
&=S^i_{\;kl}\,\varphi^k\wedge\varphi^l +R^i_{\;k\bar l}\,\varphi^k\wedge\bar\varphi^l+P^i_{\;kl}\,\varphi^k\wedge\delta v^l +Q^i_{\;k\bar l}\,\varphi^k\wedge\delta\bar v^l,
\end{aligned}
\end{equation}
where, in particular, $R^i_{\;k\bar l}=v^jR^i_{\;jk\bar l}$.

In Section~\ref{sec:3}, where considering the  special  reductive complex homogeneous spaces, we will use boldface notation, such as $\mathbf{\Gamma}^i_{\;jk}$, $\mathbf N^i_{\;k}$, $\mathbf N^i_{\;\bar k}$, $\mathbf T^i_{\;jk}$, $\bm{\omega}^i_j$, $\bm{\Omega}^i_j$ and $\bm{\Omega}^i$, to denote the corresponding quantities expressed in the adapted $(1,0)$-frame.

\section{Reductive complex Finsler homogeneous spaces}\label{sec:3}
\hspace{0.4cm}
This section is devoted to the direction-dependent complex Lie algebraic study of invariant complex Finsler metrics on reductive complex homogeneous spaces $(M=G/H,J)$.

\subsection{Invariant metrics and basic tensors}\label{sec:3.1}
\hspace{0.4cm}
Let $(M = G/H,J)$ be a reductive complex homogeneous space. For $g\in G$, we use the same symbol $g$ for the induced biholomorphism of $G/H$, given by
$$
g(aH)=gaH,\qquad a\in G,
$$
and write $dg$ for its differential. We call $(M,J,F)$ a homogeneous complex Finsler space if the action of $G$ preserves the complex Finsler metric $F$. More precisely, for $v = y - \sqrt{-1}Jy \in T^{1,0}M$ with $y \in TM$, the complex-linear extension of the real differential $dg: TM \to TM$ acts on $v$ by
\begin{equation}\label{eq:3.1}
dg(v):= dg(y)-\sqrt{-1}\,dg(Jy)= dg(y)-\sqrt{-1}\,J(dg(y)),
\end{equation}
where the second equality uses the fact that $dg$ commutes with $J$. A direct computation shows $J(dg(v))=\sqrt{-1}dg(v)$, so $dg(v)\in T^{1,0}M$. The complex Finsler metric $F$ on $M$ is $G$-invariant if for every $g \in G$,
\begin{equation}\label{eq:3.2}
F(g(z);dg(v)) = F(z;v), \qquad \forall g\in G, \forall v \in T_z^{1,0}M,\forall z\in M.
\end{equation}

A strongly pseudoconvex complex Minkowski norm $F$ on $\mathfrak m^{1,0}$ is called $\operatorname{Ad}_{\mathfrak{m}^{1,0}}(H)$-invariant, if $F$ satisfies
\begin{equation}\label{eq:3.3}
F(\operatorname{Ad}_{\mathfrak m^{1,0}}(h)u) =F(u), \qquad \forall h \in H, \; u \in \mathfrak m^{1,0}.
\end{equation}

We first have the following correspondence of invariant complex Finsler metrics on reductive complex homogeneous spaces.

\begin{proposition}[\cite{Deng2012}, \cite{LuoXiaoZhong2026}]\label{pro:3.1}
Let $(M=G/H,J)$ be a reductive complex homogeneous space. Then there is a one-to-one correspondence between $G$-invariant complex Finsler metrics on $M$  and $\operatorname{Ad}_{\mathfrak{m}^{1,0}}(H)$-invariant strongly pseudoconvex complex Minkowski norms on $\mathfrak m^{1,0}$.
\end{proposition}

\begin{remark}
Define $\mathfrak a:=\mathfrak h^{\mathbb C}\oplus\mathfrak m^{1,0}$. Using the integrability of $J$ and the $\operatorname{Ad}(H)$-invariance of $I$, it follows that $\mathfrak a$ is a complex Lie subalgebra of $\mathfrak g^{\mathbb C}$. Since $\sigma(\mathfrak a)=\mathfrak h^{\mathbb C}\oplus\mathfrak m^{0,1}$ and
$$
\mathfrak g^{\mathbb C}=\mathfrak h^{\mathbb C}\oplus\mathfrak m^{1,0}\oplus\mathfrak m^{0,1},
$$
we have $\mathfrak a\cap\sigma(\mathfrak a)=\mathfrak h^{\mathbb C}$ and $\mathfrak a+\sigma(\mathfrak a)=\mathfrak g^{\mathbb C}$.
Hence, the natural projection induces an isomorphism
$$
\mathfrak m^{1,0}\cong\mathfrak a/\mathfrak h^{\mathbb C}
$$
of complex vector spaces.  Thus Proposition~\ref{pro:3.1} is comparable to the algebraic description in \cite[Theorem~5.23]{Deng2012}.
%In the proof of \cite[Theorem~5.23]{Deng2012}, the invariant complex structure on $G/H$ is denoted by $I$, and the associated Koszul operator on $\mathfrak g$ is denoted by $J$. The latter corresponds to $\widetilde I$ in our notation: it vanishes on $\mathfrak h$ and restricts to the induced complex structure on $\mathfrak m$. Hence its complexification $J^{\mathbb C}$ vanishes on $\mathfrak h^{\mathbb C}$ and preserves $\mathfrak m^{\mathbb C}$. Its $\pm\sqrt{-1}$-eigenspaces therefore lie in $\mathfrak m^{\mathbb C}$ and coincide with $\mathfrak m^{1,0}$ and $\mathfrak m^{0,1}$, respectively.
 
\end{remark}

Let $F$ be an $\operatorname{Ad}_{\mathfrak{m}^{1,0}}(H)$-invariant complex Minkowski norm on $\mathfrak m^{1,0}$. We use the basis $\{\mathcal{E}_1,\dots,\mathcal{E}_m\}$ of $\mathfrak{m}^{1,0}$ fixed in Section~\ref{sec:2.1}. For $v=v^i\mathcal{E}_i \in \mathfrak m^{1,0}\setminus\{0\}$ and $u=u^i\mathcal{E}_i, w=w^i\mathcal{E}_i, y=y^i\mathcal{E}_i \in \mathfrak m^{1,0}$, we define the Hermitian inner product $g_v$ and the complex Cartan tensors $\mathcal C_v^+$, $\mathcal C_v^-$, respectively by
$$
g_v(u,w) := g_{i\bar j}\, u^i \bar w^j,
$$
\begin{equation*}
\mathcal C_v^+(u,w,y) := C_{i\bar j k}\, u^i \bar w^j y^k, \qquad \mathcal C_v^-(u,w,y) := C_{i\bar j \bar k}\, u^i \bar w^j \bar y^k,
\end{equation*}
where
$$
g_{i\bar j}(v) := \frac{\partial^2 F^2}{\partial v^i \partial \bar v^j}, \qquad C_{i\bar j k}(v) := \frac{\partial g_{i\bar j}}{\partial v^k}, \qquad C_{i\bar j \bar k}(v) := \frac{\partial g_{i\bar j}}{\partial \bar v^k}.
$$
We shall also need the symmetric product $\mathcal{S}_v$ and the conjugate symmetric product $\overline{\mathcal{S}}_v$ associated with $F$, which are defined, respectively by
\begin{equation}\label{eq:3.4}
\mathcal{S}_v(u,w):=g_{ij}(v)\,u^i w^j,\quad \overline{\mathcal{S}}_v(u,w):=g_{\bar i \bar j}(v)\,\bar{u}^i \bar{w}^j,
\end{equation}
where
$$
g_{ij}(v) := \frac{\partial^2 F^2}{\partial v^i \partial v^j},\quad g_{\bar i\bar j}(v) := \frac{\partial^2 F^2}{\partial \bar v^i \partial \bar v^j}.
$$

\subsection{Complex canonical connection and adapted $(1,0)$-frame fields}\label{sec:3.2}
\hspace{0.4cm}
Let $(M=G/H,J)$ be a reductive complex homogeneous space with reductive decomposition $\mathfrak g=\mathfrak h\oplus\mathfrak m$ and $\operatorname{Ad}(H)\mathfrak m\subset\mathfrak m$.

Fix the origin $o=eH$  and denote by $\mbox{exp}$ the exponential map of $G$. By Lemma 4.1 and Theorem 4.2 in \cite[Chapter II]{Helgason1978}, there exists a neighborhood $\mathcal{U}$ of $0$ in $\mathfrak m$ such that the map
$$
\Phi:\mathcal{U}\to M,\qquad \Phi(X):=\exp(X)\cdot o,
$$
is a diffeomorphism onto a neighborhood $\widetilde{\mathcal{U}}:=\Phi(\mathcal{U})$ of $o$. In particular, every point $z\in\widetilde{\mathcal{U}}$ can be uniquely written as $z=\exp(X)\cdot o$ with $X\in \mathcal{U}$. After shrinking $\mathcal U$ if necessary, we may assume that it is star-shaped with respect to $0$, so that $tX\in\mathcal U$ for every $X\in\mathcal U$ and every $t\in[0,1]$.

For $Y\in\mathfrak m$, we define an infinitesimal transformation $\iota(Y)$ on $\widetilde{\mathcal U}$ by
\begin{equation*}
\iota(Y)_z:=d(\exp(X)\circ \pi)(Y)=d\exp(X)\big\vert_o\,d\pi(Y),\qquad z=\exp(X)\cdot o\in\widetilde{\mathcal U},
\end{equation*}
where $d\pi:\mathfrak m\to T_oM$ is the isomorphism established in Section \ref{sec:2.1}. Here $d\exp(X)\big\vert_o$ denotes the differential at $o$ of the action of the group element $\exp(X)$ on $M$. At the origin $o$, we have $\iota(Y)_o=d\pi(Y).$

\begin{theorem}[\cite{Nomizu_K}, Theorem 10.2 \& Theorem 10.3]\label{th:3.3}
On a reductive homogeneous space $M$ with reductive decomposition $\mathfrak g=\mathfrak h\oplus\mathfrak m$ and $\operatorname{Ad}(H)\mathfrak m\subset\mathfrak m$, there exists a unique $G$-invariant real affine connection $\overset{\mathbf{C}}{\nabla}$ on $M$, called the canonical affine connection of the second kind, characterized by the following properties at the origin $o=eH$,
$$
\overset{\mathbf{C}}{\nabla}_{\iota(X)}\iota(Y)\big\vert_o=0,\qquad X,Y\in\mathfrak m.
$$
Its real torsion and real curvature at $o$ are given by
$$
\overset{\mathbf{C}}{T}(\iota(X),\iota(Y))_o=-\iota([X,Y]_{\mathfrak m})_o,\qquad (\overset{\mathbf{C}}{R}(\iota(X),\iota(Y))\iota(Z))_o=-\iota([[X,Y]_{\mathfrak h},Z])_o,
$$
and $\overset{\mathbf{C}}{\nabla}\, \overset{\mathbf{C}}{T}\equiv 0$, for all $X,Y,Z\in\mathfrak m$, where $[\cdot,\cdot]_{\mathfrak m}$ and $[\cdot,\cdot]_{\mathfrak h}$ denote the projections onto $\mathfrak m$ and $\mathfrak h$, respectively.
\end{theorem}

Extend $\iota$ complex-linearly to $\mathfrak m^{\mathbb C}$, still denoted by $\iota$. For $Y=Y^i(E_i-\sqrt{-1}IE_i)\in\mathfrak m^{1,0}$ with $E_i\in\mathfrak m$, we have
$$
d\pi(Y)=d\pi\big(Y^i(E_i-\sqrt{-1}IE_i)\big)=Y^id\pi(E_i)-Y^i\sqrt{-1}J_od\pi(E_i).
$$
Using the complex-linear extension of $dg$ as in \eqref{eq:3.1}, we extend $\iota(Y)$ to $\widetilde{\mathcal{U}}$ by
$$
\iota(Y)_{z}:=dg(d\pi(Y)),\quad z=gH\in\widetilde{\mathcal{U}},\quad g=\exp(X),\quad X\in\mathcal{U}.
$$
A direct computation gives $J_{z}\big(\iota(Y)_{z}\big)=\sqrt{-1}\iota(Y)_{z}.$ Hence $\iota(Y)_{z}\in T_{z}^{1,0}M$ for every $z\in\widetilde{\mathcal{U}}$. For $Y\in\mathfrak m^{1,0}$, $\overline{\iota(Y)}_{z}:=\iota(\sigma(Y))_{z}$. Then
$$
\overline{\iota(Y)}_{z}=\iota(\sigma(Y))_{z}=dg\big(d\pi(\sigma(Y))\big)=dg\big(d\pi(\overline{Y}^i(E_i+\sqrt{-1}IE_i))\big),
$$
where $\sigma$ is the conjugation of $\mathfrak g^{\mathbb C}$. In particular, $\overline{\iota(Y)}_{z}\in T_{z}^{0,1}M$.

\begin{corollary}\label{co:3.4}
Under the same conditions as in Theorem \ref{th:3.3}, if $M$, in addition, is real even-dimensional and equipped with a $G$-invariant complex structure $J$, then $\overset{\mathbf{C}}{\nabla} J\equiv0$.  Furthermore, let $\overset{\mathbf{C}}{T}$ and $\overset{\mathbf{C}}{R}$ still denote the complex-linear extensions of the real torsion and real curvature tensors of $\overset{\mathbf{C}}{\nabla}$ to $T_{\mathbb C}M$. Then the complex torsion and complex curvature tensors of the complexified canonical connection $\overset{\mathbf{C}}{\nabla}$ on $T_{\mathbb C}M$ at $o$ are given by
\begin{equation}
    \begin{split}\label{eq:3.5}
\overset{\mathbf{C}}{T}(\iota(X),\iota(Y))_o = -\iota([X,Y]_{\mathfrak m^{1,0}})_o,\qquad
\overset{\mathbf{C}}{T}(\iota(X),\overline{\iota(Y)})_o = -\iota([X,\sigma(Y)]_{\mathfrak m^{\mathbb C}})_o,
    \end{split}
\end{equation}
and
\begin{equation}
    \begin{split}\label{eq:3.6}
(\overset{\mathbf{C}}{R}(\iota(X),\iota(Y))\iota(Z))_o &= -\iota([[X,Y]_{\mathfrak h^{\mathbb C}},Z])_o,\\
(\overset{\mathbf{C}}{R}(\iota(X),\overline{\iota(Y)})\iota(Z))_o &= -\iota([[X,\sigma(Y)]_{\mathfrak h^{\mathbb C}},Z])_o,\\
(\overset{\mathbf{C}}{R}(\iota(X),\iota(Y))\overline{\iota(Z)})_o &= -\iota([[X,Y]_{\mathfrak h^{\mathbb C}},\sigma(Z)])_o,
\end{split}
\end{equation}
for all $X,Y,Z\in\mathfrak m^{1,0}$. Moreover, the complex torsion is also parallel, namely $\overset{\mathbf{C}}{\nabla}\, \overset{\mathbf{C}}{T}\equiv0$.
\end{corollary}
\begin{proof}
Since $\overset{\mathbf{C}}{\nabla}$ and $J$ are $G$-invariant, the tensor $\overset{\mathbf{C}}{\nabla} J$ is also $G$-invariant. That is, for every $g\in G$ and $X,Y\in T_oM$, we have
$$
(\overset{\mathbf{C}}{\nabla} J)_{gH}(dg(X),dg(Y))=dg((\overset{\mathbf{C}}{\nabla} J)_o(X,Y)).
$$
Therefore it suffices to verify $(\overset{\mathbf{C}}{\nabla} J)_o=0$ at the origin $o=eH$. Since $J$ is $G$-invariant, we have
$$
J_{gH}\circ dg=dg\circ J_o,\qquad \forall g\in G.
$$
Therefore, for $X\in\mathfrak m$ ,
\begin{equation*}(J\iota(X))_{gH}=J_{gH}\big(dg(d\pi(X))\big)=dg\big(J_o(d\pi(X))\big)=dg\big(d\pi(IX)\big)=\iota(IX)_{gH}.
\end{equation*}
Thus
$$
J\iota(X)=\iota(IX).
$$

Now let $X,Y\in\mathfrak m$. We have
\begin{align*}
(\overset{\mathbf{C}}{\nabla} J)_o(\iota(X),\iota(Y))&=(\overset{\mathbf{C}}{\nabla}_{\iota(X)}J)\iota(Y)\big\vert_o\\
&=\overset{\mathbf{C}}{\nabla}_{\iota(X)}(J\iota(Y))\big\vert_o-J_o\big(\overset{\mathbf{C}}{\nabla}_{\iota(X)}\iota(Y)\big\vert_o\big)\\
&=\overset{\mathbf{C}}{\nabla}_{\iota(X)}(\iota(IY))\big\vert_o-J_o\big(\overset{\mathbf{C}}{\nabla}_{\iota(X)}\iota(Y)\big\vert_o\big)\\
&=0.
\end{align*}

The formulas for the complex torsion and complex curvature now follow directly by complex-linear extension of the real formulas in Theorem \ref{th:3.3}. The identity $\overset{\mathbf{C}}{\nabla}\, \overset{\mathbf{C}}{T}=0$ is preserved under complex-linear extension, so the complex torsion is also parallel.
\end{proof}

Recall from \cite[Section 10]{Nomizu_K} that the canonical connection of the second kind on a reductive homogeneous space $M=G/H$ is also characterized by the following property:

\begin{proposition}[{\cite[Section 10, (A2)]{Nomizu_K}}]\label{pro:3.5}
On a reductive homogeneous space $M=G/H$ with decomposition $\mathfrak g=\mathfrak h\oplus\mathfrak m$, for each $X\in\mathfrak m$, there exists $r>0$ such that the curve
$$
x_t:=\exp(tX)\cdot o,\qquad 0\leq t\leq r,
$$
is a geodesic of the canonical connection of the second kind emanating from $o$. Moreover, parallel transport along $x_t$ from $o$ to $x_r$ coincides with the differential of the action of $\exp(rX)$.
\end{proposition}

Recall the real basis $\{E_1,\dots,E_m,IE_1,\dots,IE_m\}$ of $\mathfrak m$ fixed in Section~\ref{sec:2.1}, with associated complex basis $\{\mathcal E_i=E_i-\sqrt{-1}IE_i\}$ of $\mathfrak m^{1,0}$. Under the isomorphism $d\pi:\mathfrak m\to T_oM$, the real basis maps to the basis $\{d\pi(E_1),\dots,d\pi(E_m),d\pi(IE_1),\dots,d\pi(IE_m)\}$ of $T_oM$. Using the relation $J_o(d\pi(X))=d\pi(IX)$ for $X\in\mathfrak m$, we define
$$
e_i(o):=d\pi(\mathcal E_i)=d\pi(E_i)-\sqrt{-1}J_o(d\pi(E_i)).
$$
Direct computation shows $J_o(e_i(o))=\sqrt{-1}e_i(o)$. Thus $\{e_i(o)\}$ is a basis of $T_o^{1,0}M$.

For each $z\in\widetilde{\mathcal U}$, let $X\in\mathcal U$ be the unique element with $z=\exp X\cdot o$. Let $x(t)=\exp(tX)\cdot o$ be the geodesic from $o$ to $z$ with respect to the canonical connection of the second kind $\overset{\mathbf{C}}{\nabla}$. Since the canonical connection of the second kind satisfies $\overset{\mathbf{C}}{\nabla} J=0$ (Corollary \ref{co:3.4}), parallel transport along $x(t)$ preserves the $(1,0)$-type. Hence we obtain a vector, denoted by $e_i(z)\in T_z^{1,0}M$, by parallel translating $e_i(o)$ along $x(t)$. By Proposition \ref{pro:3.5}, this parallel transport coincides with the differential of the action of $g=\exp X$. Thus, we have
$$
e_i(z)=dg(e_i(o))=dg(d\pi(E_i))-\sqrt{-1}J_zdg(d\pi(E_i)),
$$
where the complex-linear extension of $dg$ is given by \eqref{eq:3.1}. Here $J_z$ denotes the value of $J$ at $z$ and we have used the $G$-invariance of $J$, i.e. $dg\circ J_o=J_z\circ dg$.

Consequently, $\{e_i\}$ is a smooth local $(1,0)$-frame field on $\widetilde{\mathcal U}$, obtained by parallel translating a basis of $T_o^{1,0}M$ along geodesics. The neighborhood $\widetilde{\mathcal U}$ is a normal neighborhood of $o$ with respect to the canonical connection of the second kind: the exponential map $X\mapsto\exp X\cdot o$ is a diffeomorphism from $\mathcal U\subset\mathfrak m$ onto $\widetilde{\mathcal U}$, and each $z=\exp X\cdot o$ is joined to $o$ by the unique geodesic $x(t)=\exp(tX)\cdot o$, $0\leq t\leq 1$ (see \cite[Chapter~I, \S\,6]{Helgason1978}). Moreover, for any $z=\exp X\cdot o\in\widetilde{\mathcal{U}}$, $g=\exp X$,
$$
J_z(e_i(z))=J_z(dg(e_i(o)))=dg(J_o(e_i(o)))=dg(\sqrt{-1}e_i(o))=\sqrt{-1}e_i(z).
$$
Hence $\{e_i\}$ is indeed a smooth $(1,0)$-frame field on $\widetilde{\mathcal U}$. This adapted frame is not assumed to be holomorphic or to be a coordinate frame. Similarly, we define
$$
\bar e_i(o):=d\pi(\sigma(\mathcal E_i))=d\pi(E_i)+\sqrt{-1}J_o(d\pi(E_i))\in T_o^{0,1}M,
$$
and for $z=\exp X\cdot o\in\widetilde{\mathcal U}$,
$$
\bar e_i(z):=dg(\bar e_i(o))=dg(d\pi(E_i))+\sqrt{-1}J_z\,dg(d\pi(E_i))\in T_z^{0,1}M.
$$

\begin{proposition}\label{pro:3.6}
For a fixed vector $v_o:=v^i e_i(o)\in T_o^{1,0}M$ with $v^i\in\mathbb{C}$ and each $z=\exp X\cdot o\in\widetilde{\mathcal U}$, $X\in\mathfrak m$, define the vector $v_z:=dg(v_o)$, where $g=\exp X$. Then
$$
v_z=v^i e_i(z)\in T_z^{1,0}M.
$$
\end{proposition}
\begin{proof}
Since $v_o=v^ie_i(o)$ and $e_i(z)=dg(e_i(o))$, the natural complex-linear extension of $dg$ gives
$$
v_z=dg(v_o)=v^i dg(e_i(o))=v^ie_i(z).
$$
Moreover, $J_z(v_z)=\sqrt{-1}v_z$ because each $e_i(z)$ is of type $(1,0)$. Hence $v_z\in T_z^{1,0}M$ for every $z\in\widetilde{\mathcal{U}}$.
\end{proof}

\subsection{Complex Rund connection and the structure functions}\label{sec:3.3}
\hspace{0.4cm}
Now we investigate the relationship between the coefficients of the complex Rund connection and the structure constants of $\mathfrak g^{\mathbb C}$. Let $(M=G/H,J)$ be a reductive complex homogeneous space and $F:T^{1,0}M\to[0,+\infty)$ a $G$-invariant complex Finsler metric on $M$. On the normal neighborhood $\widetilde{\mathcal U}$, we have the adapted (1,0) frame $\{e_i\}$ constructed in Section \ref{sec:3.2}.

By Proposition \ref{pro:3.6}, the same fiber coordinates $(v^1,\dots,v^m)$ determine vectors at every point of $\widetilde{\mathcal U}$ via the adapted (1,0) frame. The $G$-invariance of $F$ then yields, for each $z\in \widetilde{\mathcal U}$ and every $v_z\in T_z^{1,0}M$,
\begin{equation*}
F(z;v_z)=F(z;v^i e_i(z))=F(g(o);dg(v^i e_i(o)))=F(o;v^i e_i(o)).
\end{equation*}
Consequently, when expressed with respect to the adapted $(1,0)$ frame $\{e_i\}$ on $\widetilde{\mathcal U}$, the complex Finsler metric $F$ is independent of the base point $z$. In particular, the fundamental tensor $g_{i\bar j}$ and the Cartan tensors $C_{i\bar j k}$, $C_{i\bar j\bar k}$ are completely determined by their values at the origin $o$:
$$
g_{i\bar j}(z;v)=g_{i\bar j}(o;v),\qquad C_{i\bar j k}(z;v)=C_{i\bar j k}(o;v),\qquad C_{i\bar j\bar k}(z;v)=C_{i\bar j\bar k}(o;v).
$$
Equivalently, since $F$ is identified with a complex Minkowski norm on $\mathfrak m^{1,0}$ via $d\pi$, the fundamental tensor $g_{i\bar j}$ and the Cartan tensors $C_{i\bar j k}$, $C_{i\bar j\bar k}$ can be computed directly as the fiber derivatives of this norm.

\begin{remark}
Via the isomorphism $d\pi:\mathfrak m^{1,0}\to T_o^{1,0}M$ and the adapted $(1,0)$-frame $\{e_i\}$, a nonzero vector $v=v^i\mathcal E_i\in\mathfrak m^{1,0}$ is naturally identified with the vector field
$$
v^ie_i(z)\in T_z^{1,0}M,\qquad z\in\widetilde{\mathcal U},
$$
whose value at $o$ is $d\pi(v)=v^ie_i(o)$. Thus, in the complex homogeneous case, the same symbol $v$ stands for both the vector in $\mathfrak m^{1,0}$ and its associated local $(1,0)$-vector field on $\widetilde{\mathcal U}$. We call $v$ the \emph{reference vector} at which the fiber derivatives of $F$ are evaluated.
\end{remark}

The complex Rund connection is defined on the pulled-back bundle $\tilde{\pi}^\ast T^{1,0}M$. Pulling back the local smooth adapted $(1,0)$ frame $\{e_i\}$ to $\tilde{\pi}^\ast T^{1,0}M$ gives a local smooth frame for $\tilde{\pi}^\ast T^{1,0}M$, still denoted by $\{e_i\}$, with dual $(1,0)$-forms $\{\varphi^i\}$ satisfying $\varphi^i(e_j)=\delta^i_j$. The connection 1-forms $\{\bm{\omega}^i_j\}$ are uniquely determined by the first structure equation
\begin{equation}\label{eq:3.7}
d\varphi^i=\varphi^j\wedge\bm{\omega}^i_j+\bm{\tau}^i,
\end{equation}
where $\bm{\tau}^i$ is the complex torsion form of the complex Rund connection and contains only horizontal $(2,0)$-components, and the almost metric compatibility
\begin{equation}\label{eq:3.8}
dg_{i\bar j}=g_{l\bar j}\,\bm{\omega}^l_i+g_{i\bar l}\,\bm{\omega}^{\bar l}_{\bar j}+C_{i\bar j k}\,\delta v^k+C_{i\bar j\bar k}\,\delta\bar v^k,
\end{equation}
where $\delta v^i=dv^i+v^j\bm{\omega}^i_j$.

For $z\in\widetilde{\mathcal U}$, let $\lambda^i_{\;jk}(z)$, $\lambda^i_{\;j\bar k}(z)$, and $\lambda^{\bar i}_{\;j\bar k}(z)$ be the complex structure functions defined by the Lie brackets of the adapted (1,0) frame $\{e_i\}$ and its conjugate $\{\bar e_i\}$:
$$
[e_j,e_k]=\lambda^i_{\;jk}\,e_i,\qquad [e_j,\bar e_k]=\lambda^i_{\;j\bar k}\,e_i+\lambda^{\bar i}_{\;j\bar k}\,\bar e_i.
$$
Equivalently, for the dual coframe $\{\varphi^i\}$,
\begin{equation}\label{eq:3.9}
d\varphi^i=-\frac12\lambda^i_{\;jk}\,\varphi^j\wedge\varphi^k-\lambda^i_{\;j\bar k}\,\varphi^j\wedge\bar\varphi^k.
\end{equation}

Write $\bm{\omega}^i_{j} = \mathbf{\Gamma}^i_{\;jk}\varphi^k + \mathbf{\Gamma}^i_{\;j\bar k}\bar\varphi^k$ and $\bm{\tau}^i=\mathbf{T}^i_{\;jk}\varphi^j\wedge\varphi^k$. Solving equations \eqref{eq:3.7}, \eqref{eq:3.8}, and \eqref{eq:3.9} simultaneously for $\mathbf{\Gamma}^i_{\;jk}$, $\mathbf{\Gamma}^i_{\;j\bar k}$, and $\mathbf{T}^i_{\;jk}$, we obtain (see \cite[Section 4.2]{LuoXiaoZhong2025} for details of the computation)
\begin{align*}
\mathbf{\Gamma}^i_{\;jk} &=\lambda^{\;\,i}_{j\;k}-C^{\;i}_{j\;l}\,\mathbf{N}^l_{\;k}+C^{\;i}_{j\;\bar l}\,\lambda^{\bar l}_{\;\bar o k}, \qquad \mathbf{\Gamma}^i_{\;j\bar k}=-\lambda^i_{\;j\bar k}; \\
\mathbf{T}^i_{\;jk} &=\frac12\Bigl\{\lambda^{\;\,i}_{k\;j}-\lambda^i_{\;jk}-\lambda^{\;\,i}_{j\;k}+C^{\;i}_{j\;l}\,\mathbf{N}^l_{\;k}-C^{\;i}_{j\;\bar l}\,\lambda^{\bar l}_{\;\bar o k}-C^{\;i}_{k\;l}\,\mathbf{N}^l_{\;j}+C^{\;i}_{k\;\bar l}\,\lambda^{\bar l}_{\;\bar o j}\Bigr\},
\end{align*}
where
\begin{equation}\label{eq:3.10}
\mathbf{N}^j_{\;k}=\lambda^{\;\,j}_{o\;k}+g^{\bar q j}g_{\bar q\bar l}\lambda^{\bar l}_{\;\bar o k}.
\end{equation}
Here $\mathbf T^i_{\;jk}=-\mathbf T^i_{\;kj}$, so that $\bm{\tau}^i(e_j,e_k)=2\mathbf T^i_{\;jk}$. In particular, $2\mathbf{T}^i_{\;jk}=\mathbf{\Gamma}^i_{\;kj}-\mathbf{\Gamma}^i_{\;jk}-\lambda^i_{\;jk}$.

In the adapted frame, the forms $\delta v^i$ and the contracted complex Rund curvature forms $\bm\Omega^i$ introduced in Section~\ref{sec:2} satisfy
$$
\delta v^i=dv^i+v^j\bm\omega^i_j=dv^i+\mathbf N^i_{\;k}\varphi^k+\mathbf N^i_{\;\bar k}\bar\varphi^k,\qquad \bm\Omega^i=d(\delta v^i)-\delta v^j\wedge\bm\omega^i_j=v^j\bm\Omega^i_j.
$$
With respect to the coframe $\{\varphi^i,\bar\varphi^i,\delta v^i,\delta\bar v^i\}$ associated with the adapted frame $\{e_i,\bar{e}_i\}$ and the complex Rund connection, we compute the coefficient of $\varphi^k\wedge\bar\varphi^j$ in the contracted complex Rund curvature form $\bm\Omega^i$. Following Section~4.2 in \cite{LuoXiaoZhong2025} and using $\mathbf N^i_{\;\bar j}=-\lambda^i_{\;p\bar j}v^p$, we obtain
\begin{align*}
\mathbf{R}^i_{\;k\bar{j}}:=&e_k(\mathbf{N}^i_{\;\bar{j}})-\bar{e}_j(\mathbf{N}^i_{\;k})-[\mathbf{N}^i_{\;k}]_{v^l}\,\lambda^l_{\;p\bar{j}}v^p+[\mathbf{N}^i_{\;k}]_{\bar v^l}\,\mathbf{N}^{\bar l}_{\;\bar{j}}-\mathbf{N}^i_{\;l}\,\lambda^l_{\;k\bar{j}}+\lambda^i_{\;l\bar{j}}\,\mathbf{N}^l_{\;k}\nonumber\\
&-\lambda^i_{\;p\bar l}v^p\,\lambda^{\bar l}_{\;\bar{j} k},
\end{align*}
where $\mathbf{N}^{\bar l}_{\;\bar{j}}=\overline{\mathbf{N}^{l}_{\;j}}$. Here $[\cdot]_{v^l}$ and $[\cdot]_{\bar v^l}$ denote fiber derivatives, while $e_k(\cdot)$ and $\bar e_j(\cdot)$ denote directional derivatives in the base variable along $e_k$ and $\bar e_j$, respectively, with the adapted-frame fiber coordinates $v^i,\bar v^i$ held fixed. Compared with the Lie group case, the additional contribution consists of the directional derivative terms
\begin{equation}\label{eq:3.11}
e_k(\mathbf N^i_{\;\bar j})-\bar e_j(\mathbf N^i_{\;k}).
\end{equation}
These terms vanish in a left-invariant frame, whose structure functions are constant. For the adapted frame on $G/H$, however, the structure functions generally vary with the base point. Evaluating the terms in \eqref{eq:3.11} at $o$ therefore requires the values of the first derivatives of the structure functions along the adapted frame fields $e_k$ and $\bar e_j$ at $o$.

We decompose the contracted horizontal complex $(1,1)$-curvature coefficients $\mathbf{R}^i_{\;k\bar{j}}$ into two parts
\begin{equation}\label{eq:3.12}
\mathbf{R}^i_{\;k\bar{j}}=\mathfrak{R}^i_{\;k\bar{j}}+\mathcal{R}^i_{\;k\bar{j}},
\end{equation}
where
\begin{align}
\mathfrak{R}^i_{\;k\bar{j}}&:=e_k(\mathbf{N}^i_{\;\bar{j}})-\bar{e}_j(\mathbf{N}^i_{\;k}),\label{eq:3.13}\\
\mathcal{R}^i_{\;k\bar{j}}&:=-[\mathbf{N}^i_{\;k}]_{v^l}\,\lambda^l_{\;p\bar{j}}v^p+[\mathbf{N}^i_{\;k}]_{\bar v^l}\,\mathbf{N}^{\bar l}_{\;\bar{j}}-\mathbf{N}^i_{\;l}\,\lambda^l_{\;k\bar{j}}+\lambda^i_{\;l\bar{j}}\,\mathbf{N}^l_{\;k}-\lambda^i_{\;p\bar l}v^p\,\lambda^{\bar l}_{\;\bar{j} k}.\label{eq:3.14}
\end{align}
The technical lemma in the next subsection expresses these data in terms of the complex torsion $T$, complex curvature $R$ and covariant derivative $\nabla T$ at $o$ of the complexified $J$-preserving connection used to construct the adapted frame. We then apply these identities to the canonical connection of the second kind on a reductive complex homogeneous space, whose complex torsion is parallel and whose complex torsion and complex curvature are determined by Lie brackets, to obtain the required Lie algebra formulas.

\subsection{A technical lemma on the adapted (1,0) frame}
\hspace{0.4cm}
Let $M$ be an even-dimensional smooth real manifold equipped with a complex structure $J$, and let $\nabla$ be a real affine connection on $TM$ satisfying $\nabla J=0$. We regard $\nabla$ as a complex connection on $T_{\mathbb C}M$ preserving the decomposition
$$
T_{\mathbb C}M=T^{1,0}M\oplus T^{0,1}M.
$$
Denote its complex torsion and complex curvature tensors by $T$ and $R$, respectively. They are the complex-linear extensions of the real torsion and real curvature tensors of the original affine connection on $TM$.

Let $\{e_i\}$ be the local smooth $(1,0)$-frame on a normal neighborhood of $o\in M$, obtained by parallel transporting a basis of $T_o^{1,0}M$ along the radial geodesics of $\nabla$.

Define functions $\lambda^i_{\;jk}$, $\lambda^i_{\;j\bar k}$, $\lambda^{\bar i}_{\;j\bar k}$ satisfying
$$
[e_j,e_k]=\lambda^i_{\;jk}e_i,\qquad [e_j,\bar e_k]=\lambda^i_{\;j\bar k}e_i+\lambda^{\bar i}_{\;j\bar k}\bar e_i.
$$

We use the following component notation for $T$, $R$ and $\nabla T$:
\begin{align}
    &T(e_j,e_k) =T^i_{\;jk}e_i,\quad T(e_j,\bar{e}_k) =T^i_{\;j\bar{k}}e_i+T^{\bar{i}}_{\;j\bar{k}}\bar{e}_i,\quad (\nabla_{e_l}T)(e_j,e_k)=\nabla_{l}T^i_{\;jk}e_i,\nonumber\\
    &R(e_j,e_k)e_l=R^i_{\;ljk}e_i,\quad R(\bar{e}_j,e_k)e_l=R^i_{\;l\bar{j}k}e_i,\quad R(\bar{e}_j,\bar{e}_k)e_l=R^i_{\;l\bar{j}\bar{k}}e_i,\label{eq:3.15}\\
    &(\nabla_{e_l}T)(e_j,\bar{e}_k)=\nabla_{l}T^i_{\;j\bar k}e_i+\nabla_{l}T^{\bar{i}}_{\;j\bar k}\bar{e}_i,\quad (\nabla_{\bar{e}_l}T)(e_j,\bar{e}_k)=\nabla_{\bar{l}}T^i_{\;j\bar k}e_i+\nabla_{\bar{l}}T^{\bar{i}}_{\;j\bar k}\bar{e}_i.\nonumber
\end{align}

\begin{lemma}\label{lm:3.9}
With the notation above, the following identities hold at $o$:
$$
\lambda^i_{\;jk}=-T^i_{\;jk},\qquad \lambda^i_{\;j\bar k}=-T^i_{\;j\bar k},\qquad \lambda^{\bar i}_{\;j\bar k}=-T^{\bar i}_{\;j\bar k},
$$
and
\begin{align*}
e_l(\lambda^i_{\;jk})&=\frac12(R^i_{\;klj}-R^i_{\;jlk})-\nabla_l T^i_{\;jk},& \bar e_l(\lambda^i_{\;jk})&=\frac12(R^i_{\;k\bar lj}-R^i_{\;j\bar lk})-\nabla_{\bar l}T^i_{\;jk},\\
e_l(\lambda^i_{\;j\bar k})&=\frac12 R^i_{\;j\bar kl}-\nabla_l T^i_{\;j\bar k},& \bar e_l(\lambda^i_{\;j\bar k})&=\frac12 R^i_{\;j\bar k\bar l}-\nabla_{\bar l}T^i_{\;j\bar k}.
\end{align*}
\end{lemma}

\begin{remark}
By Corollary~\ref{co:3.4}, the canonical connection $\overset{\mathbf{C}}{\nabla}$ of the second kind satisfies $\overset{\mathbf{C}}{\nabla}J=0$. Hence Lemma~\ref{lm:3.9} applies to reductive complex homogeneous spaces equipped with this connection.
\end{remark}

\begin{remark}
The proof combines the radial parallel transport identity with complex polarization and the complex curvature relations for $\nabla$ to determine the first derivatives of the connection coefficients at $o$. The values and first derivatives of the adapted-frame structure functions then follow from the complex torsion identity and its covariant derivative.
\end{remark}

\begin{proof}
Choose a real basis $\{E_1,\dots,E_m,J_oE_1,\dots,J_oE_m\}$ of $T_oM$ such that
\begin{equation}\label{eq:3.16}
e_i(o)=E_i-\sqrt{-1}J_oE_i,\qquad
\bar e_i(o)=E_i+\sqrt{-1}J_oE_i.
\end{equation}
For a constant parameter $\zeta=(\zeta^1,\dots,\zeta^m)\in\mathbb C^m$, define the real initial velocity
$$
X_\zeta:=\zeta^i e_i(o)+\bar\zeta^i\bar e_i(o)\in T_oM.
$$
By \eqref{eq:3.16}, every real initial velocity has this form. Let $\operatorname{Exp}^{\nabla}_o$ denote the exponential map of $\nabla$. After shrinking the normal neighborhood, for every $\zeta$ in a sufficiently small ball centered at $0\in\mathbb C^m$, the radial geodesic
$$
\gamma_\zeta(t):=\operatorname{Exp}^{\nabla}_o(tX_\zeta),\qquad 0\leq t\leq1,
$$
is defined and remains in this neighborhood. The numbers $\zeta^i$ are parameters of the initial velocity, not coordinate functions on $M$.

The fields $e_i$ are obtained by parallel transporting $e_i(o)$ along these radial geodesics. Since $\nabla J=0$, parallel transport preserves the type. Thus $\{e_i,\bar e_i\}$ is a smooth frame of $T_{\mathbb C}M$ of types $(1,0)$ and $(0,1)$; no holomorphicity of the frame is assumed.

Let $\{\varphi^i,\bar\varphi^i\}$ be the coframe dual to $\{e_i,\bar e_i\}$ on $T_{\mathbb C}M$. Thus
$$
\varphi^i(e_j)=\delta^i_j,\quad \varphi^i(\bar e_j)=0,\qquad \bar\varphi^i(e_j)=0,\quad \bar\varphi^i(\bar e_j)=\delta^i_j.
$$
Define the complex-valued connection $1$-forms $\omega^i_j$ and the complex-valued coefficient functions $\Gamma^i_{\;jk},\Gamma^i_{\;j\bar k}$, respectively, by
\begin{equation}\label{eq:3.17}
\begin{aligned}
\nabla_Xe_j=\omega^i_j(X)e_i,\qquad \omega^i_j=\Gamma^i_{\;jk}\varphi^k+\Gamma^i_{\;j\bar k}\bar\varphi^k,\quad X\in T_{\mathbb{C}}M
\end{aligned}
\end{equation}

The velocity $\dot\gamma_\zeta$ and the field
$\zeta^ke_k(\gamma_\zeta(t))+\bar\zeta^k\bar e_k(\gamma_\zeta(t))$
are parallel along $\gamma_\zeta$ and have the same initial value $X_\zeta$. Uniqueness for parallel transport therefore gives
\begin{equation}\label{eq:3.18}
\dot\gamma_\zeta(t)=\zeta^ke_k(\gamma_\zeta(t))+\bar\zeta^k\bar e_k(\gamma_\zeta(t)).
\end{equation}
Consequently, $\varphi^k(\dot\gamma_\zeta(t))=\zeta^k$ and $\bar\varphi^k(\dot\gamma_\zeta(t))=\bar\zeta^k$. Using $\nabla_{\dot\gamma_\zeta}e_j=0$ in \eqref{eq:3.17}, we obtain
\begin{equation}\label{eq:3.19}
\zeta^k\Gamma^i_{\;jk}(\gamma_\zeta(t))+\bar\zeta^k\Gamma^i_{\;j\bar k}(\gamma_\zeta(t))=0.
\end{equation}
At $t=0$, equation \eqref{eq:3.19} reduces to
$$
\zeta^k\Gamma^i_{\;jk}(o)+\bar\zeta^k\Gamma^i_{\;j\bar k}(o)=0.
$$
The parameter $\sqrt{-1}\zeta$ determines another real radial geodesic $\gamma_{\sqrt{-1}\zeta}$, whose initial velocity is $X_{\sqrt{-1}\zeta}=J_oX_\zeta$. Applying the same radial parallel transport identity to this curve at $t=0$ and dividing by $\sqrt{-1}$ gives
$$
\zeta^k\Gamma^i_{\;jk}(o)-\bar\zeta^k\Gamma^i_{\;j\bar k}(o)=0.
$$
Using the arbitrariness of $\zeta$, we conclude that
\begin{equation}\label{eq:3.20}
\Gamma^i_{\;jk}(o)=\Gamma^i_{\;j\bar k}(o)=0,\qquad \omega^i_j\big\vert_o=0.
\end{equation}
Thus, for every complex vector field $X$, we have
$$
\nabla_Xe_j\big\vert_o=0,\qquad \nabla_X\bar e_j\big\vert_o=0.
$$

For every smooth complex-valued function $f$, \eqref{eq:3.18} implies
$$
\frac{d}{dt}f(\gamma_\zeta(t))
=\zeta^le_l(f)(\gamma_\zeta(t))+\bar\zeta^l\bar e_l(f)(\gamma_\zeta(t)).
$$
Differentiate \eqref{eq:3.19} with respect to $t$ and set $t=0$. This yields
\begin{equation}\label{eq:3.21}
\begin{aligned}
0={}\zeta^k\zeta^l e_l(\Gamma^i_{\;jk})(o)+\zeta^k\bar\zeta^l\bar e_l(\Gamma^i_{\;jk})(o)+\bar\zeta^k\zeta^l e_l(\Gamma^i_{\;j\bar k})(o)+\bar\zeta^k\bar\zeta^l\bar e_l(\Gamma^i_{\;j\bar k})(o).
\end{aligned}
\end{equation}
By homogeneity, this quadratic identity holds for every $\zeta\in\mathbb C^m$.

Fix $i,j$. Introduce the following abbreviations for complex-valued quadratic polynomials in the parameter $\zeta$:
\begin{align*}
A^i_j(\zeta):&=\zeta^k\zeta^l e_l(\Gamma^i_{\;jk})(o),\quad B^i_j(\zeta):=\bar \zeta^k\bar \zeta^l\bar e_l(\Gamma^i_{\;j\bar k})(o),\\
W^i_j(\zeta):&=\zeta^k\bar \zeta^l\bigl\{\bar e_l(\Gamma^i_{\;jk})(o)+e_k(\Gamma^i_{\;j\bar l})(o)\bigr\}.
\end{align*}
Relabeling the summation indices $k,l$ in the third term of \eqref{eq:3.21}, that equation reduces to $A^i_j(\zeta)+W^i_j(\zeta)+B^i_j(\zeta)=0$.

Replacing $\zeta$ by $\sqrt{-1}\zeta$ in $A^i_j(\zeta)+W^i_j(\zeta)+B^i_j(\zeta)=0$ gives a second identity
\begin{equation*}
-A^i_j(\zeta)+W^i_j(\zeta)-B^i_j(\zeta)=0,
\end{equation*}
because $A^i_j(\sqrt{-1}\zeta)=-A^i_j(\zeta)$, $W^i_j(\sqrt{-1}\zeta)=W^i_j(\zeta)$, and $B^i_j(\sqrt{-1}\zeta)=-B^i_j(\zeta)$.
Thus, we have $W^i_j(\zeta)=0$ and $A^i_j(\zeta)+B^i_j(\zeta)=0$. The parameter $e^{\sqrt{-1}\pi/4}\zeta$ likewise determines a real radial geodesic $\gamma_{e^{\sqrt{-1}\pi/4}\zeta}$. Replacing $\zeta$ by $e^{\sqrt{-1}\pi/4}\zeta$ in $A^i_j(\zeta)+B^i_j(\zeta)=0$ gives
$$
\sqrt{-1}A^i_j(\zeta)-\sqrt{-1}B^i_j(\zeta)=0.
$$
Together with $A^i_j(\zeta)+B^i_j(\zeta)=0$, this proves
\begin{equation*}
A^i_j(\zeta)=W^i_j(\zeta)=B^i_j(\zeta)=0.
\end{equation*}

Let $\epsilon_k$ denote the $k$-th standard basis vector of $\mathbb C^m$. In each test, the selected indices $k,l$ are fixed and are not summed.

From $A^i_j(\zeta)=0$, taking $\zeta=\epsilon_k$ gives $e_k(\Gamma^i_{\;jk})(o)=0$. For $k\neq l$, taking $\zeta=\epsilon_k+\epsilon_l$ and removing the two diagonal terms gives
\begin{equation}\label{eq:3.22}
e_l(\Gamma^i_{\;jk})(o)+e_k(\Gamma^i_{\;jl})(o)=0.
\end{equation}
The same identity holds for $k=l$ by the diagonal calculation. Likewise, applying these two tests to $B^i_j(\zeta)=0$ gives
\begin{equation}\label{eq:3.23}
\bar e_l(\Gamma^i_{\;j\bar k})(o)+\bar e_k(\Gamma^i_{\;j\bar l})(o)=0
\end{equation}
for all $k,l$.

For the mixed polynomial, abbreviate
$$
(W^i_j)_{k\bar l}:=\bar e_l(\Gamma^i_{\;jk})(o)+e_k(\Gamma^i_{\;j\bar l})(o).
$$
The equality $W^i_j(\zeta)=0$ reads $\zeta^k\bar \zeta^l (W^i_j)_{k\bar l}=0$. Taking $\zeta=\epsilon_k$ gives $(W^i_j)_{k\bar k}=0$. For $k\neq l$, the tests $\zeta=\epsilon_k+\epsilon_l$ and $\zeta=\epsilon_k+\sqrt{-1}\epsilon_l$ give, respectively,
$$
(W^i_j)_{k\bar l}+(W^i_j)_{l\bar k}=0,\qquad-\sqrt{-1} (W^i_j)_{k\bar l}+\sqrt{-1} (W^i_j)_{l\bar k}=0.
$$
Hence $(W^i_j)_{k\bar l}=(W^i_j)_{l\bar k}=0$. Consequently,
\begin{equation}\label{eq:3.24}
\bar e_l(\Gamma^i_{\;jk})(o)+e_k(\Gamma^i_{\;j\bar l})(o)=0.
\end{equation}
Interchanging $k,l$ also gives the useful equivalent form
\begin{equation}\label{eq:3.25}
e_l(\Gamma^i_{\;j\bar k})(o)+\bar e_k(\Gamma^i_{\;jl})(o)=0.
\end{equation}

For arbitrary smooth complex vector fields $X,Y$, the definition of the connection forms gives
\begin{align*}
\nabla_X\nabla_Ye_j&=\nabla_X\bigl(\omega^p_j(Y)e_p\bigr)=\bigl\{X(\omega^i_j(Y))+\omega^p_j(Y)\omega^i_p(X)\bigr\}e_i,\\
\nabla_Y\nabla_Xe_j&=\bigl\{Y(\omega^i_j(X))+\omega^p_j(X)\omega^i_p(Y)\bigr\}e_i,\\
\nabla_{[X,Y]}e_j&=\omega^i_j([X,Y])e_i.
\end{align*}
It follows that
\begin{align*}
R(X,Y)e_j=\bigl\{X(\omega^i_j(Y))-Y(\omega^i_j(X))+\omega^p_j(Y)\omega^i_p(X)-\omega^p_j(X)\omega^i_p(Y)-\omega^i_j([X,Y])\bigr\}e_i.
\end{align*}
By \eqref{eq:3.20}, all product terms and the last term of the above equation vanish at $o$. Therefore
\begin{equation}\label{eq:3.26}
R(X,Y)e_j\big\vert_o=\bigl\{X(\omega^i_j(Y))-Y(\omega^i_j(X))\bigr\}(o)e_i(o).
\end{equation}

With the component convention \eqref{eq:3.15}, applying \eqref{eq:3.26} to the pairs $(X,Y)=(e_l,e_k)$ and $(X,Y)=(\bar e_l,e_k)$ gives
\begin{align*}
R^i_{\;jlk}(o)&=e_l(\Gamma^i_{\;jk})(o)-e_k(\Gamma^i_{\;jl})(o)=2e_l(\Gamma^i_{\;jk})(o),\\
R^i_{\;j\bar l k}(o)&=\bar e_l(\Gamma^i_{\;jk})(o)-e_k(\Gamma^i_{\;j\bar l})(o)=2\bar e_l(\Gamma^i_{\;jk})(o),
\end{align*}
where the last equalities use \eqref{eq:3.22} and \eqref{eq:3.24}, respectively. Applying the same calculation to $(e_l,\bar e_k)$ and $(\bar e_l,\bar e_k)$ gives
\begin{align*}
R^i_{\;jl\bar k}(o)&=e_l(\Gamma^i_{\;j\bar k})(o)-\bar e_k(\Gamma^i_{\;jl})(o)=2e_l(\Gamma^i_{\;j\bar k})(o),\\
R^i_{\;j\bar l\bar k}(o)&=\bar e_l(\Gamma^i_{\;j\bar k})(o)-\bar e_k(\Gamma^i_{\;j\bar l})(o)=2\bar e_l(\Gamma^i_{\;j\bar k})(o),
\end{align*}
where the last equalities use  \eqref{eq:3.25} and \eqref{eq:3.23}. Thus all four derivative identities are
\begin{equation}\label{eq:3.27}
\begin{aligned}
e_l(\Gamma^i_{\;jk})(o)&=\frac12R^i_{\;jlk}(o),\quad
\bar e_l(\Gamma^i_{\;jk})(o)=\frac12R^i_{\;j\bar l k}(o),\\
e_l(\Gamma^i_{\;j\bar k})(o)&=\frac12R^i_{\;jl\bar k}(o),\quad
\bar e_l(\Gamma^i_{\;j\bar k})(o)=\frac12R^i_{\;j\bar l\bar k}(o).
\end{aligned}
\end{equation}

In the frame $\{e_i,\bar e_i\}$, the defining relations for the structure functions and for the components of the complex torsion tensor $T$ of $\nabla$ are as follows:
\begin{align*}
[e_j,e_k]&=\lambda^i_{\;jk}e_i,\quad [e_j,\bar e_k]=\lambda^i_{\;j\bar k}e_i+\lambda^{\bar i}_{\;j\bar k}\bar e_i,\\
T(e_j,e_k)&=T^i_{\;jk}e_i,\quad T(e_j,\bar e_k)=T^i_{\;j\bar k}e_i+T^{\bar i}_{\;j\bar k}\bar e_i.
\end{align*}
At $o$, all frame derivatives under $\nabla$ vanish by \eqref{eq:3.20} and its conjugate. The complex torsion identity therefore gives
$$
T(e_j,e_k)\big\vert_o=-[e_j,e_k]\big\vert_o,\qquad T(e_j,\bar e_k)\big\vert_o=-[e_j,\bar e_k]\big\vert_o.
$$
Comparing the $(1,0)$- and $(0,1)$-components in these two identities yields
\begin{align*}
\lambda^i_{\;jk}(o)=-T^i_{\;jk}(o),\quad \lambda^i_{\;j\bar k}(o)=-T^i_{\;j\bar k}(o),\quad \lambda^{\bar i}_{\;j\bar k}(o)=-T^{\bar i}_{\;j\bar k}(o).
\end{align*}

Expanding $(\nabla_{e_l}T)(e_j,e_k)$ gives
\begin{align*}
(\nabla_{e_l}T)(e_j,e_k)={}&\nabla_{e_l}\bigl(T(e_j,e_k)\bigr)-T(\nabla_{e_l}e_j,e_k)-T(e_j,\nabla_{e_l}e_k).
\end{align*}
The last two terms of $(\nabla_{e_l}T)(e_j,e_k)$ vanish at $o$ and the first reduces to
\begin{align*}
\nabla_{e_l}(T^i_{\;jk}e_i)\big\vert_o=e_l(T^i_{\;jk})(o)e_i(o)+T^i_{\;jk}(o)\nabla_{e_l}e_i\big\vert_o=e_l(T^i_{\;jk})(o)e_i(o).
\end{align*}
For the mixed inputs, the same expansion has both $e_i$ and $\bar e_i$ output terms, whose connection contributions also vanish at $o$. Replacing the derivative direction by $\bar e_l$ has the same effect. Hence
\begin{equation}\label{eq:3.28}
\begin{aligned}
e_l(T^i_{\;jk})(o)&=\nabla_lT^i_{\;jk}(o),\quad \bar e_l(T^i_{\;jk})(o)=\nabla_{\bar l}T^i_{\;jk}(o),\\
e_l(T^i_{\;j\bar k})(o)&=\nabla_lT^i_{\;j\bar k}(o),\quad \bar e_l(T^i_{\;j\bar k})(o)=\nabla_{\bar l}T^i_{\;j\bar k}(o).
\end{aligned}
\end{equation}

Throughout the neighborhood, the complex $(2,0)$-torsion identity is
\begin{align*}
T(e_j,e_k)=\nabla_{e_j}e_k-\nabla_{e_k}e_j-[e_j,e_k]=\bigl(\Gamma^i_{\;kj}-\Gamma^i_{\;jk}-\lambda^i_{\;jk}\bigr)e_i.
\end{align*}
Thus
\begin{equation}\label{eq:3.29}
\lambda^i_{\;jk}=\Gamma^i_{\;kj}-\Gamma^i_{\;jk}-T^i_{\;jk}.
\end{equation}
Differentiating \eqref{eq:3.29} along $e_l$ and using \eqref{eq:3.27} and \eqref{eq:3.28}, we obtain
\begin{align*}
e_l(\lambda^i_{\;jk})(o)=e_l(\Gamma^i_{\;kj})(o)-e_l(\Gamma^i_{\;jk})(o)-e_l(T^i_{\;jk})(o)=\frac12R^i_{\;klj}(o)-\frac12R^i_{\;jlk}(o)-\nabla_lT^i_{\;jk}(o).
\end{align*}
Differentiating \eqref{eq:3.29} along $\bar e_l$ instead gives
\begin{align*}
\bar e_l(\lambda^i_{\;jk})(o)=\bar e_l(\Gamma^i_{\;kj})(o)-\bar e_l(\Gamma^i_{\;jk})(o)-\bar e_l(T^i_{\;jk})(o)=\frac12R^i_{\;k\bar l j}(o)-\frac12R^i_{\;j\bar l k}(o)-\nabla_{\bar l}T^i_{\;jk}(o).
\end{align*}

Taking the $(1,0)$-component of $T(e_j,\bar e_k)=\nabla_{e_j}\bar e_k-\nabla_{\bar e_k}e_j-[e_j,\bar e_k]$ therefore yields $T^i_{\;j\bar k}=-\Gamma^i_{\;j\bar k} -\lambda^i_{\;j\bar k}$, or equivalently
\begin{equation*}
\lambda^i_{\;j\bar k}=-\Gamma^i_{\;j\bar k}-T^i_{\;j\bar k}.
\end{equation*}
Differentiate the above equation along $e_l$ and evaluate at $o$:
\begin{align*}
e_l(\lambda^i_{\;j\bar k})(o)=-e_l(\Gamma^i_{\;j\bar k})(o)-e_l(T^i_{\;j\bar k})(o)=\frac12R^i_{\;j\bar k l}(o)-\nabla_lT^i_{\;j\bar k}(o).
\end{align*}
The last equality follows from
$R(e_l,\bar e_k)=-R(\bar e_k,e_l)$. Similarly,
\begin{align*}
\bar e_l(\lambda^i_{\;j\bar k})(o)=-\bar e_l(\Gamma^i_{\;j\bar k})(o)-\bar e_l(T^i_{\;j\bar k})(o)=\frac12R^i_{\;j\bar k\bar l}(o)-\nabla_{\bar l}T^i_{\;j\bar k}(o),
\end{align*}
where we used $R(\bar e_l,\bar e_k)=-R(\bar e_k,\bar e_l)$.

This completes the proof.
\end{proof}

\begin{remark}
For the corresponding real structure functions and their derivatives in the real tangent bundle $TM$, we refer the reader to \cite[Lemma 4.4]{Huang2015} for detailed computations and explicit formulas.
\end{remark}

\subsection{K\"ahler--Finsler metrics and weakly K\"ahler--Finsler metrics}
\hspace{0.4cm}
We now apply Lemma \ref{lm:3.9} to the reductive complex homogeneous space $(M=G/H,J)$ equipped with the canonical connection of the second kind.

Let $\{e_i\}$ be the adapted $(1,0)$-frame constructed in Section \ref{sec:3.2}. Let $\lambda^i_{\;jk}$, $\lambda^i_{\;j\bar k}$, $\lambda^{\bar i}_{\;j\bar k}$ be the complex structure functions of $\{e_i,\bar{e}_i\}$ defined by
$$
[e_j,e_k]=\lambda^i_{\;jk}e_i,\qquad [e_j,\bar e_k]=\lambda^i_{\;j\bar k}e_i+\lambda^{\bar i}_{\;j\bar k}\bar e_i.
$$
The constants $c^i_{\;jk}$, $c^i_{\;j\bar k}$, $c^{\bar i}_{\;j\bar k}$, $c^\alpha_{\;l\bar j}$, $c^k_{\;\alpha i}$ are the structure constants of $\mathfrak g^{\mathbb C}$ defined in Section \ref{sec:2.1}.

\begin{lemma}\label{lm:3.13}
With the notation above, at the origin $o=eH$, the values of these structure functions and their derivatives along the frame directions are given by
\begin{align}
\lambda^i_{\;jk}(o)&=c^i_{\;jk},\qquad\lambda^i_{\;j\bar k}(o)=c^i_{\;j\bar k},\qquad\lambda^{\bar i}_{\;j\bar k}(o)=c^{\bar i}_{\;j\bar k},\label{eq:3.30}\\
e_l(\lambda^i_{\;j\bar k})(o)&=-\frac12 c^\alpha_{\;\bar k l}\,c^i_{\;\alpha j},\qquad \bar e_l(\lambda^i_{\;j\bar k})(o)=-\frac12 c^\alpha_{\;\bar k \bar l}\,c^i_{\;\alpha j},\label{eq:3.31}\\
e_l(\lambda^i_{\;jk})(o)&=-\frac12 (c^\alpha_{\;lj}\,c^i_{\;\alpha k}-c^\alpha_{\;lk}\,c^i_{\;\alpha j}),\quad \bar e_l(\lambda^i_{\;jk})(o)=-\frac12 (c^\alpha_{\;\bar{l}j}\,c^i_{\;\alpha k}-c^\alpha_{\;\bar{l}k}\,c^i_{\;\alpha j}).\label{eq:3.32}
\end{align}
\end{lemma}

\begin{proof}
By Corollary \ref{co:3.4}, the canonical connection $\overset{\mathbf{C}}{\nabla}$ of the second kind is complex-linear, and its complex torsion and complex curvature tensors at $o$ are expressed purely in terms of the Lie algebra data via \eqref{eq:3.5} and \eqref{eq:3.6}. Explicitly, we have
\begin{align*}
  \overset{\mathbf{C}}{T}\,{}^i_{\;jk}(o)=-c^i_{\;jk},\quad \overset{\mathbf{C}}{T}\,{}^i_{\;j\bar{k}}(o)&=-c^i_{\;j\bar{k}},\quad \overset{\mathbf{C}}{T}\,{}^{\bar{i}}_{\;j\bar{k}}(o)=-c^{\bar{i}}_{\;j\bar{k}},\\
  \overset{\mathbf{C}}{R}\,{}^i_{\;jlk}(o)=-c^{\alpha}_{\;lk}\,c^i_{\;\alpha j},\quad \overset{\mathbf{C}}{R}\,{}^i_{\;j\bar lk}(o)&=-c^{\alpha}_{\;\bar{l}k}\,c^i_{\;\alpha j},\quad \overset{\mathbf{C}}{R}\,{}^i_{\;j\bar l\bar k}(o)=-c^{\alpha}_{\;\bar{l}\bar{k}}\,c^i_{\;\alpha j}, \\
  \overset{\mathbf{C}}{\nabla}\, \overset{\mathbf{C}}{T}&\equiv0.
\end{align*}
Substituting the above equations into Lemma \ref{lm:3.9} gives \eqref{eq:3.30}, \eqref{eq:3.31} and \eqref{eq:3.32}.
\end{proof}

Now, substituting \eqref{eq:3.30} into the general formulas for the complex Rund connection coefficients $\mathbf{\Gamma}^i_{\;jk}, \mathbf{\Gamma}^i_{\;j\bar k}$ and the horizontal complex $(2,0)$-torsion coefficients $\mathbf{T}^i_{\;jk}$ (obtained in Section \ref{sec:3.3}) yields their explicit values at the origin. In the adapted frame we obtain
\begin{align}
\mathbf{\Gamma}^i_{\;jk}(o;v)&=c^{\;\,i}_{j\;k}-C^{\;i}_{j\;l}\,\mathbf{N}^l_{\;k}+C^{\;i}_{j\;\bar l}\,c^{\bar l}_{\;\bar o k},\qquad \mathbf{\Gamma}^i_{\;j\bar k}(o;v)=-c^i_{\;j\bar k};\label{eq:3.33}\\
\mathbf{T}^i_{\;jk}(o;v)&=\frac12\bigl\{c^{\;\,i}_{k\;j}-c^i_{\;jk}-c^{\;\,i}_{j\;k}+C^{\;i}_{j\;l}\,\mathbf{N}^l_{\;k}-C^{\;i}_{j\;\bar l}\,c^{\bar l}_{\;\bar o k}-C^{\;i}_{k\;l}\,\mathbf{N}^l_{\;j}+C^{\;i}_{k\;\bar l}\,c^{\bar l}_{\;\bar o j}\bigr\},\label{eq:3.34}
\end{align}
where
\begin{equation}
\mathbf{N}^j_{\;k}(o;v)=c^{\;\,j}_{o\;k}+g^{\bar q j}g_{\bar q\bar l}\,c^{\bar l}_{\;\bar o k}.\label{eq:3.35}
\end{equation}

Recall from Proposition \ref{pro:3.1} that the $G$-invariant complex Finsler metric $F$ on $M$ corresponds to an $\operatorname{Ad}_{\mathfrak{m}^{1,0}}(H)$-invariant complex Minkowski norm on $\mathfrak m^{1,0}$ via the isomorphism $d\pi:\mathfrak m^{1,0}\to T_o^{1,0}M$. Under this identification, we may regard the fundamental tensor $g_{i\bar j}$, the complex Cartan tensors $C_{i\bar{j}l}$, $C_{i\bar{j}\bar{l}}$, and the connection coefficients $\mathbf{N}^j_{\;k}$ in \eqref{eq:3.33}--\eqref{eq:3.35} as quantities defined on $\mathfrak m^{1,0}$. Consequently, all subsequent computations will be carried out on $\mathfrak m^{1,0}$.
\begin{definition}\label{df:3.14}
The connection operator $\mathcal N_v:\mathfrak m^{1,0}\to\mathfrak m^{1,0}$ is defined by
$$
\mathcal N_v(w):= \mathbf{N}_{\;j}^i(v) w^j \mathcal{E}_i,\quad \forall w=w^i\mathcal{E}_i\in\mathfrak{m}^{1,0},
$$
where $v\in\mathfrak m^{1,0}\setminus\{0\}$ is the reference vector and the coefficients $\mathbf{N}_{\;j}^i(v)$ are given by \eqref{eq:3.35}.
\end{definition}
\begin{proposition}
The operator $\mathcal N_v:\mathfrak m^{1,0}\to\mathfrak m^{1,0}$ defined in Definition \ref{df:3.14} satisfies the invariant identity
\begin{equation}\label{eq:3.36}
g_v(\mathcal N_v(w),u) = g_v(v,[u,\sigma(w)]_{\mathfrak{m}^{1,0}}) + \overline{\mathcal{S}}_v(u,[v,\sigma(w)]_{\mathfrak{m}^{1,0}})
\end{equation}
for all $u,w\in\mathfrak m^{1,0}$.
\end{proposition}
\begin{proof}
    Expanding \eqref{eq:3.35} gives
\begin{equation*}
    \mathbf{N}^j_{\;k}=g^{\bar{p}j}g_{l\bar{i}}\,c^{\bar{i}}_{\;\bar{p}k}v^l+g^{\bar p j}g_{\bar p\bar l}\,c^{\bar l}_{\;\bar i k}\bar{v}^i.
\end{equation*}
Multiplying the preceding identity by $g_{j\bar q}$ and summing over $j$, we obtain
\begin{equation*}
    g_{j\bar{q}}\mathbf{N}^j_{\;k}=g_{l\bar{i}}\,c^{\bar{i}}_{\;\bar{q}k}v^l+g_{\bar q\bar l}\,c^{\bar l}_{\;\bar{i} k}\bar{v}^i.
\end{equation*}
Contracting the above identity with $w^k$ and $\bar{u}^q$ and using the definitions of $g_v(\cdot,\cdot)$ and $\overline{\mathcal{S}}_v(\cdot,\cdot)$, we obtain \eqref{eq:3.36}.
\end{proof}

Thus, by the definitions of a K\"ahler--Finsler metric and a weakly K\"ahler--Finsler metric, we have the following theorem.

\begin{theorem}\label{th:3.16}
Let $(M=G/H,J)$ be a reductive complex homogeneous space with reductive decomposition $\mathfrak g=\mathfrak h\oplus\mathfrak m$ and $\operatorname{Ad}(H)\mathfrak m\subset\mathfrak m$. Let $F:T^{1,0}M\to[0,\infty)$ be a $G$-invariant complex Finsler metric on $M$. Then
\begin{enumerate}
    \item[(1)] $F$ is a K\"ahler--Finsler metric if and only if
$$
c_{k\bar q o}-c_{\bar q ok}-c_{o\bar q k}-g_{\bar q\bar l}\,c^{\bar l}_{\ \bar o k}-C_{k\bar q l}\,\mathbf{N}^l_{\;o}+C_{k\bar q\bar l}\,c^{\bar l}_{\ \bar o o}=0,
$$
or equivalently, for all $u,w\in\mathfrak m^{1,0}$ and $v\in\mathfrak m^{1,0}\setminus\{0\}$,
\begin{align}
&g_v(w,[u,\sigma(v)]_{\mathfrak m^{1,0}})-g_v([v,w]_{\mathfrak m^{1,0}},u)-g_v(v,[u,\sigma(w)]_{\mathfrak m^{1,0}})\nonumber\\
&-\overline{\mathcal{S}}_v(u,[v,\sigma(w)]_{\mathfrak m^{1,0}})-\mathcal C_v^+(w,u,\mathcal N_v(v))+\mathcal C_v^-(w,u,[v,\sigma(v)]_{\mathfrak m^{1,0}})=0;\label{eq:3.37}
\end{align}
    \item[(2)] $F$ is a weakly K\"ahler--Finsler metric if and only if
$$
c_{k\bar o o}-c_{\bar o ok}-c_{o\bar o k}-g_{kl}\,\mathbf{N}^l_{\;o}=0,
$$
or equivalently, for all $w\in\mathfrak m^{1,0}$ and $v\in\mathfrak m^{1,0}\setminus\{0\}$,
\begin{align}
g_v(w,[v,\sigma(v)]_{\mathfrak m^{1,0}})-g_v([v,w]_{\mathfrak m^{1,0}},v)-g_v(v,[v,\sigma(w)]_{\mathfrak m^{1,0}})-\mathcal{S}_v(\mathcal N_v(v),w)=0.\label{eq:3.38}
\end{align}
\end{enumerate}
\end{theorem}

\begin{proof}
We identify $T_o^{1,0}M$ with $\mathfrak m^{1,0}$ via $d\pi$ and evaluate all coefficients at $(o;v)$, where $v\in\mathfrak m^{1,0}\setminus\{0\}$ is arbitrary.

Multiplying \eqref{eq:3.34} by $g_{i\bar q}v^j$ and using $v^j C_{j\bar q l}=0$, we obtain the identity
\begin{equation}
2g_{i\bar q}\,\mathbf T^i_{\;jk}v^j=\Bigl\{g_{k\bar l}\,c^{\bar l}_{\;\bar q j}-g_{l\bar q}\,c^l_{\;jk}-g_{j\bar l}\,c^{\bar l}_{\;\bar q k}-C_{j\bar q\bar l}\bar v^p \,c^{\bar l}_{\;\bar p k}-C_{k\bar q l}\,\mathbf N^l_{\;j}+C_{k\bar q\bar l}\bar v^p\, c^{\bar l}_{\;\bar p j}\Bigr\}v^j.\label{eq:3.39}
\end{equation}

By the result of Chen and Shen \cite{ChenAndShen}, the K\"ahler--Finsler condition $\mathbf T^i_{\;jk}=0$ is equivalent to $\mathbf T^i_{\;jk}v^j=0$ for every nonzero $v$. Since $(g_{i\bar q})$ is positive definite, this is equivalent to the vanishing of the right-hand side of \eqref{eq:3.39}. Contract it with $w^k\bar u^q$, where $u,w\in\mathfrak m^{1,0}$ are arbitrary. The first three terms give the first three inner products in \eqref{eq:3.37}, using the conjugation relations for the structure constants. For the remaining terms, $v^j C_{j\bar q\bar l}=g_{\bar q\bar l}$ and Definition~\ref{df:3.14} give
\begin{gather*}
g_{\bar q\bar l}\,\bar u^q\bar v^p\, c^{\bar l}_{\;\bar p k}w^k=\overline{\mathcal S}_v\bigl(u,[v,\sigma(w)]_{\mathfrak m^{1,0}}\bigr),\\
C_{k\bar q l}\,w^k\bar u^q\,\mathbf N^l_{\;i}v^i=\mathcal C_v^+\bigl(w,u,\mathcal N_v(v)\bigr),\quad C_{k\bar q\bar l}\,w^k\bar u^q\bar v^p \,c^{\bar l}_{\;\bar p i}v^i=\mathcal C_v^-\bigl(w,u,[v,\sigma(v)]_{\mathfrak m^{1,0}}\bigr).
\end{gather*}
Thus the resulting scalar identity is exactly \eqref{eq:3.37}. Conversely, if it holds for all $u,w$ and every nonzero $v$, the right-hand side of \eqref{eq:3.39} vanishes componentwise, proving the K\"ahler--Finsler characterization.

For the weakly K\"ahler--Finsler condition, contract the identity \eqref{eq:3.39} with $\bar v^q$. The relation $\bar v^q C_{j\bar q\bar l}=0$ eliminates both terms involving $C_{j\bar q\bar l}$ or $C_{k\bar q\bar l}$, while $\bar v^q C_{k\bar q l}=g_{kl}$ gives
\begin{align*}
2g_{i\bar q}\,\mathbf T^i_{\;jk}v^j\bar v^q=\bigl\{g_{k\bar l}\,c^{\bar l}_{\;\bar q j}-g_{l\bar q}\,c^l_{\;jk}-g_{j\bar l}\,c^{\bar l}_{\;\bar q k}\bigr\}v^j\bar v^q-g_{kl}\,\mathbf N^l_{\;j}v^j.
\end{align*}
By definition, $F$ is weakly K\"ahler--Finsler if and only if the left-hand side vanishes for every $k$ and every nonzero $v$. Contracting with arbitrary $w^k$, the first three terms reduce to the corresponding inner products in \eqref{eq:3.38}. Since $g_{kl}=g_{lk}$, the last term gives
$$
-g_{kl}\,w^k\mathbf N^l_{\;i}v^i=-\mathcal S_v\bigl(\mathcal N_v(v),w\bigr).
$$
This proves \eqref{eq:3.38}, and the converse follows from the arbitrariness of $w$.
\end{proof}

\subsection{Holomorphic curvatures}
\hspace{0.4cm}
We now explain how the curvature on $\tilde\pi^*T^{1,0}M$ is expressed on $\mathfrak m^{1,0}$. Since $F$ is $G$-invariant, each $g\in G$ induces a bundle map, for $(z;v)\in\tilde{M}$ and $\xi\in T^{1,0}_M$
$$
((z;v),\xi)\longmapsto((g(z);dg_zv),dg_z\xi)
$$
covering $(z;v)\mapsto(g(z);dg_zv)$ on $\tilde M$. By the characterization of the complex Rund connection, this map preserves the connection, its horizontal lifts and its curvature. Transitivity therefore reduces the curvature calculation to all points $(o;v)$ with $v\in T_o^{1,0}M\setminus\{0\}$.

For the quotient map $\pi:G\to G/H$, the complex-linear extension of $d\pi$ restricts to an isomorphism $\mathfrak m^{1,0}\to T_o^{1,0}M$. Hence
\begin{align*}
\tilde\pi^*T^{1,0}M\Big|_{T_o^{1,0}M\setminus\{0\}}\cong\bigl(T_o^{1,0}M\setminus\{0\}\bigr)\times T_o^{1,0}M\cong\bigl(\mathfrak m^{1,0}\setminus\{0\}\bigr)\times\mathfrak m^{1,0}.
\end{align*}
Explicitly, for $v\in\mathfrak m^{1,0}\setminus\{0\}$ and $u\in\mathfrak m^{1,0}$, this identification is given by
$$
\bigl((o;d\pi(v)),d\pi(u)\bigr)\longleftrightarrow (v,u).
$$
In what follows, we use the same symbols $v,u,w$ for their images under $d\pi$, and write $(o;v)$ for the corresponding point of $\tilde M$.

For $w\in\mathfrak m^{1,0}$, the horizontal lifts $\chi_{(o;v)}(w)$ and $\chi_{(o;v)}(\sigma(w))$ belong to $(\mathcal H_{\mathbb C})_{(o;v)}$. We evaluate the curvature on these lifts and on $u$ in the pullback fiber:
$$
\Omega_{(o;v)} \bigl(\chi_{(o;v)}(w),\chi_{(o;v)}(\sigma(w))\bigr)u\in \bigl(\tilde\pi^*T^{1,0}M\bigr)_{(o;v)} \cong T_o^{1,0}M \cong\mathfrak m^{1,0}.
$$

Recall that $\bm\Omega^i=v^l\bm\Omega^i_l$ by \eqref{eq:2.2}. Let $\mathbf R^i_{\;lk\bar j}$ be the coefficient of $\varphi^k\wedge\bar\varphi^j$ in $\bm\Omega^i_l$ with respect to the associated coframe. For each reference vector $v\in\mathfrak m^{1,0}\setminus\{0\}$, the horizontal complex $(1,1)$-curvature operator is
\begin{align*}
\mathbf R_v(w,\sigma(w))u
&:=\Omega_{(o;v)}\bigl(\chi_{(o;v)}(w),\chi_{(o;v)}(\sigma(w))\bigr)u\\
&=\mathbf R^i_{\;lk\bar j}(o;v)\,u^lw^k\bar w^j\mathcal E_i, \qquad u,w\in\mathfrak m^{1,0},
\end{align*}
where the first equality is understood under the identifications above. In particular,
$$
\mathbf R_v(w,\sigma(w))v =\mathbf R^i_{\;lk\bar j}(o;v)\,v^lw^k\bar w^j\mathcal E_i=\mathbf R^i_{\;k\bar j}(o;v)\,w^k\bar w^j\mathcal E_i.
$$
By \eqref{eq:3.12}, this contracted curvature decomposes as
$$
\mathbf R_v(w,\sigma(w))v=\mathfrak R_v(w,\sigma(w))v+\mathcal R_v(w,\sigma(w))v,
$$
where
\begin{align*}
\mathfrak R_v(w,\sigma(w))v:=\mathfrak R^i_{\;k\bar j}(o;v)\,w^k\bar w^j\mathcal E_i,\quad\mathcal R_v(w,\sigma(w))v:=\mathcal R^i_{\;k\bar j}(o;v)\,w^k\bar w^j\mathcal E_i.
\end{align*}
The two contributions are defined by \eqref{eq:3.13} and \eqref{eq:3.14} in the chosen adapted frame.

By Lemma~\ref{lm:3.13}, the contracted horizontal curvature coefficients at $(o;v)$ are determined by the structure constants of $\mathfrak g^{\mathbb C}$ and the invariant complex Minkowski norm on $\mathfrak m^{1,0}$. In particular, substituting \eqref{eq:3.30} into \eqref{eq:3.14} gives
\begin{align*}
\mathcal{R}^i_{\;k\bar{j}}&=-[\mathbf{N}^i_{\;k}]_{v^l}\,c^l_{\;p\bar{j}}v^p+[\mathbf{N}^i_{\;k}]_{\bar v^p}\,\mathbf{N}^{\bar p}_{\;\bar{j}}-\mathbf{N}^i_{\;l}\,c^l_{\;k\bar{j}}+c^i_{\;l\bar{j}}\,\mathbf{N}^l_{\;k}-c^i_{\;p\bar l}v^p\,c^{\bar l}_{\;\bar{j} k}.
\end{align*}

We now express the fiber derivatives in the curvature formula as directional derivatives on $\mathfrak m^{1,0}\setminus\{0\}$. Throughout these derivatives, the base point is fixed at $o$, while the reference vector $v$ varies. On the trivial bundle $(\mathfrak m^{1,0}\setminus\{0\})\times\mathfrak m^{1,0}$, let $\mathcal D$ be the standard flat connection induced by the vector space $\mathfrak m^{1,0}$. Let $\{\mathcal E_i\}$ be a basis of $\mathfrak m^{1,0}$ and write $v=v^i\mathcal E_i$, so that $(v^1,\dots,v^m)$ are global holomorphic coordinates and $\{\partial_{v^k}:=\frac{\partial}{\partial v^k}\}$ is a holomorphic coordinate frame. Under the natural identifications of the holomorphic and antiholomorphic tangent spaces with $\mathfrak m^{1,0}$ and $\mathfrak m^{0,1}$, respectively, $w=w^k\mathcal E_k$ and $\sigma(w)$ correspond to $w^k\partial_{v^k}$ and $\bar w^k\partial_{\bar v^k}$. The connection $\mathcal D$, extended complex-linearly in the direction argument, makes the constant frame $\{\mathcal E_i\}$ parallel in both types of directions. For a smooth complex-valued function $f(v,\bar v)$ and a smooth $\mathfrak m^{1,0}$-valued function $\phi (v,\bar v)=\phi^i(v,\bar v)\mathcal E_i$, we therefore have
\begin{equation}\label{eq:3.40}
\begin{aligned}
\mathcal D_w f&=w^k\frac{\partial f}{\partial v^k},&\mathcal D_{\sigma(w)}f&=\bar w^k\frac{\partial f}{\partial\bar v^k},\\
\mathcal D_w \phi&=w^k\frac{\partial \phi^i}{\partial v^k}\mathcal E_i,&\mathcal D_{\sigma(w)}\phi&=\bar w^k\frac{\partial \phi^i}{\partial\bar v^k}\mathcal E_i.
\end{aligned}
\end{equation}
The same formulas apply pointwise when $w$ is a smooth direction field.

As the reference vector varies, $v=v^i\mathcal E_i$ defines a section of $(\mathfrak m^{1,0}\setminus\{0\})\times\mathfrak m^{1,0}$, given by $v\mapsto(v,v)$. We continue to denote this section by $v$ and call it the
\emph{radial vector with variable direction}. By \eqref{eq:3.40},
$$
\mathcal D_w v=w,\qquad \mathcal D_{\sigma(w)}v=0.
$$
When differentiating $v\mapsto\mathcal N_v(w)$, the argument $w$ is held fixed. For $v\mapsto\mathcal N_v(v)$, the product rule also differentiates the last argument $v$.

A computation using the connection operator $\mathcal N_v$ and the directional derivative $\mathcal D$ yields
\begin{align}
\mathcal{R}_v&(w, \sigma(w))v\nonumber\\
=&w^k\bar w^j\left\{-[\mathbf{N}^i_{\;k}]_{v^l}\,c^l_{\;p \bar{j}}v^p+[\mathbf{N}^i_{\;k}]_{\bar v^l}\,\mathbf{N}^{\bar l}_{\;\bar{j}}-\mathbf{N}^i_{\;l}\,c^l_{\;k \bar{j}}+c^i_{\;l \bar{j}}\,\mathbf{N}^l_{\;k}-c^i_{\;p\bar l}v^p\,c^{\bar l}_{\;\bar{j} k}\right\}\mathcal{E}_i\nonumber\\
=&\mathcal{D}_{\sigma(\mathcal N_v(w))}\mathcal N_v(w)-\mathcal{D}_{[v, \sigma(w)]_{\mathfrak{m}^{1, 0}}}\mathcal N_v(w)-\mathcal N_v([w, \sigma(w)]_{\mathfrak{m}^{1, 0}})\nonumber\\
&+[\mathcal N_v(w), \sigma(w)]_{\mathfrak{m}^{1, 0}}-[v, [ \sigma(w), w]_{\mathfrak{m}^{0,1}}]_{\mathfrak{m}^{1, 0}}.\label{eq:3.41}
\end{align}
Setting $w=v$ and using the product rule
$$
\frac{\partial}{\partial v^i}\bigl(\mathbf N^j_{\;p}(v)v^p\bigr)=\frac{\partial\mathbf N^j_{\;p}}{\partial v^i}(v)v^p+\mathbf N^j_{\;i}(v),
$$
one obtains
\begin{equation}
\begin{split}
\mathcal{R}_v(v, \sigma(v))v=&\mathcal{D}_{\sigma(\mathcal N_v(v))}\mathcal N_v(v)-\mathcal{D}_{[v, \sigma(v)]_{\mathfrak{m}^{1, 0}}}\mathcal N_v(v)+[\mathcal N_v(v), \sigma (v)]_{\mathfrak{m}^{1, 0}}\\
&-[v, [ \sigma(v), v]_{\mathfrak{m}^{0,1}}]_{\mathfrak{m}^{1, 0}}.\label{eq:3.42}
\end{split}
\end{equation}

The curvature operator $\mathfrak R$ corresponding to the derivative contribution is independent of the invariant complex Finsler metric, as the following proposition shows.

\begin{proposition}\label{pro:3.17}
Let $(M=G/H,J)$ be a reductive complex homogeneous space equipped with a $G$-invariant strongly pseudoconvex complex Finsler metric $F$. For every reference vector $v\in\mathfrak m^{1,0}\setminus\{0\}$ and all $u,w,y\in\mathfrak m^{1,0}$, the curvature operator $\mathfrak R_v$ corresponding to the derivative contribution satisfies
\begin{equation}\label{eq:3.43}
\mathfrak R_v(u,\sigma(w))y =[[\sigma(w),u]_{\mathfrak h^{\mathbb C}},y]_{\mathfrak m^{1,0}}.
\end{equation}
\end{proposition}

\begin{proof}
In the chosen adapted frame, we first introduce the uncontracted derivative coefficients
\begin{equation}\label{eq:3.44}
\mathfrak R^i_{\;lk\bar j}(o;v):=e_k(\mathbf\Gamma^i_{\;l\bar j})(o;v)-\bar e_j(\mathbf\Gamma^i_{\;lk})(o;v).
\end{equation}
Put
$$
\mathfrak R_v(u,\sigma(w))y:=\mathfrak R^i_{\;lk\bar j}(o;v)\,y^l u^k\bar w^j\mathcal E_i,\qquad u,w,y\in\mathfrak m^{1,0}.
$$
Since $\mathbf N^i_{\;k}=v^l\mathbf\Gamma^i_{\;lk}$ and $\mathbf N^i_{\;\bar j}=v^l\mathbf\Gamma^i_{\;l\bar j}$, we have $v^l\mathfrak R^i_{\;lk\bar j}=\mathfrak R^i_{\;k\bar j}$, so this definition agrees with \eqref{eq:3.13} upon contraction with $v$.

Recall that $\{\mathcal F_\alpha\}$ is a basis of $\mathfrak h$. For each fixed $\alpha$, the $\operatorname{Ad}_{\mathfrak m^{1,0}}(H)$-invariance of $F$ gives
$$
F^2\left(\operatorname{Ad}_{\mathfrak m^{1,0}}\bigl(\exp(t\mathcal F_\alpha)\bigr)v\right)=F^2(v),\qquad t\in\mathbb R.
$$
Differentiating the above invariance identity at $t=0$ gives
\begin{equation}\label{eq:3.45}
\frac{\partial F^2}{\partial v^p}\,c^p_{\;\alpha q}v^q+\frac{\partial F^2}{\partial\bar v^l}\,c^{\bar l}_{\;\alpha\bar q}\bar v^q=0,
\end{equation}
where we used the following identities
\begin{align*}
\frac{d}{dt}\Big|_{t=0}\left(\operatorname{Ad}_{\mathfrak m^{1,0}}\bigl(\exp(t\mathcal F_\alpha)\bigr)v\right)^p=c^p_{\;\alpha l}v^l,\quad \frac{d}{dt}\Big|_{t=0}\overline{\left(\operatorname{Ad}_{\mathfrak m^{1,0}}\bigl(\exp(t\mathcal F_\alpha)\bigr)v\right)^q}=c^{\bar q}_{\;\alpha\bar l}\bar v^l.
\end{align*}

Now differentiating \eqref{eq:3.45} with respect to $\bar v^q$ and using $\frac{\partial F^2}{\partial\bar v^l}=g_{p\bar l}v^p$, we obtain
\begin{equation}\label{eq:3.46}
g_{p\bar q}\,c^p_{\;\alpha j}v^j+g_{p\bar k}v^p\,c^{\bar k}_{\;\alpha\bar q}+g_{\bar q\bar k}\,c^{\bar k}_{\;\alpha\bar j}\bar v^j=0.
\end{equation}
Multiplying on both sides of the above equation by $g^{\bar q i}$ gives
\begin{equation}\label{eq:3.47}
g^{\bar q i}\Bigl(g_{p\bar l}v^p\,c^{\bar l}_{\;\alpha\bar q}+g_{\bar q\bar l}\,c^{\bar l}_{\;\alpha\bar p}\bar v^p\Bigr)=-c^i_{\;\alpha p}v^p.
\end{equation}

To obtain the identity for an arbitrary last argument, we differentiate \eqref{eq:3.46} once more with respect to $v^l$ and use $C_{p\bar q l}v^p=0$ and obtain
\begin{align*}
0=g_{p\bar q}\,c^p_{\;\alpha l}+g_{l\bar p}\,c^{\bar p}_{\;\alpha\bar q}+C_{l\bar q p}\,c^p_{\;\alpha j}v^j+C_{l\bar q\bar p}\,c^{\bar p}_{\;\alpha\bar j}\bar v^j.
\end{align*}
Multiplying on both sides of the above equation by $g^{\bar q i}$ using the raised Cartan components gives
\begin{equation}\label{eq:3.48}
\begin{aligned}
g^{\bar q i}g_{l\bar p}\,c^{\bar p}_{\;\alpha\bar q}+C^{\;i}_{l\;p}\,c^p_{\;\alpha j}v^j+C^{\;i}_{l\;\bar p}\,c^{\bar p}_{\;\alpha\bar j}\bar v^j=-c^i_{\;\alpha l}.
\end{aligned}
\end{equation}

By \eqref{eq:3.31} and taking its complex conjugate, we have
\begin{equation}
e_k(\lambda^i_{\;l\bar j})(o)=-\frac12c^\alpha_{\;\bar jk} c^i_{\;\alpha l},\quad \bar{e}_k(\lambda^{\bar{i}}_{\;\bar{l}j})(o)=-\frac12c^\alpha_{\;j\bar{k}} c^{\bar{i}}_{\;\alpha \bar{l}}.\label{eq:3.49}
\end{equation}

Expanding \eqref{eq:3.10} gives
\begin{equation}\label{eq:3.50}
\mathbf N^i_{\;k}=g^{\bar q i}\Bigl(g_{p\bar l}v^p\lambda^{\bar l}_{\;\bar q k}+g_{\bar q\bar l}\lambda^{\bar l}_{\;\bar p k}\bar v^p\Bigr).
\end{equation}
Differentiating \eqref{eq:3.50} in the direction $\bar e_j$ and applying the second equality of  \eqref{eq:3.49} gives
\begin{align}
\bar e_j(\mathbf N^i_{\;k})(o;v)&=g^{\bar q i}\Bigl(g_{p\bar l}v^p\bar e_j(\lambda^{\bar l}_{\;\bar q k})(o)+g_{\bar q\bar l}\bar e_j(\lambda^{\bar l}_{\;\bar p k})(o)\bar v^p\Bigr)\nonumber\\
&=\frac12c^\alpha_{\;\bar jk} g^{\bar q i}\Bigl(g_{p\bar l}v^p c^{\bar l}_{\;\alpha\bar q}+g_{\bar q\bar l}c^{\bar l}_{\;\alpha\bar p}\bar v^p\Bigr)\nonumber\\
&=-\frac12c^\alpha_{\;\bar jk} c^i_{\;\alpha p}v^p,\label{eq:3.51}
\end{align}
where the last equality follows from \eqref{eq:3.47}.

The coefficients obtained in Section~\ref{sec:3.3}, using the same raised Cartan components, are
\begin{equation}\label{eq:3.52}
\begin{aligned}
\mathbf\Gamma^i_{\;lk}=g^{\bar q i}g_{l\bar p}\lambda^{\bar p}_{\;\bar q k}-C^{\;i}_{l\;p}\mathbf N^p_{\;k}+C^{\;i}_{l\;\bar p}\lambda^{\bar p}_{\;\bar q k}\bar v^q,\quad\mathbf\Gamma^i_{\;l\bar j}&=-\lambda^i_{\;l\bar j}.
\end{aligned}
\end{equation}
The first term of \eqref{eq:3.44} follows immediately from the first equality of \eqref{eq:3.49}:
\begin{equation}\label{eq:3.53}
e_k(\mathbf\Gamma^i_{\;l\bar j})(o;v)=\frac12c^\alpha_{\;\bar jk}c^i_{\;\alpha l}.
\end{equation}
Differentiating \eqref{eq:3.52} along $\bar{e}_j$ and using the second equality of  \eqref{eq:3.49} and \eqref{eq:3.51} gives
\begin{align}
\bar e_j(\mathbf\Gamma^i_{\;lk})(o;v)={}&g^{\bar q i}g_{l\bar p}\,\bar e_j(\lambda^{\bar p}_{\;\bar q k})(o)-C^{\;i}_{l\;p}\,\bar e_j(\mathbf N^p_{\;k})(o;v)+C^{\;i}_{l\;\bar p}\,\bar e_j(\lambda^{\bar p}_{\;\bar q k})(o)\bar v^q\nonumber\\
={}&\frac12c^\alpha_{\;\bar jk}g^{\bar q i}g_{l\bar p}c^{\bar p}_{\;\alpha\bar q}
+\frac12c^\alpha_{\;\bar jk}C^{\;i}_{l\;p}c^p_{\;\alpha q}v^q+\frac12c^\alpha_{\;\bar jk}C^{\;i}_{l\;\bar p}c^{\bar p}_{\;\alpha\bar q}\bar v^q\nonumber\\
={}&\frac12c^\alpha_{\;\bar jk}\Bigl(g^{\bar q i}g_{l\bar p}c^{\bar p}_{\;\alpha\bar q}+C^{\;i}_{l\;p}c^p_{\;\alpha q}v^q+C^{\;i}_{l\;\bar p}c^{\bar p}_{\;\alpha\bar q}\bar v^q\Bigr)\nonumber\\
={}&-\frac12c^\alpha_{\;\bar jk}c^i_{\;\alpha l},\label{eq:3.54}
\end{align}
where the last equality follows from \eqref{eq:3.48}.

Equations \eqref{eq:3.53} and \eqref{eq:3.54} express the required quantities entirely in terms of the structure constants of $\mathfrak g^{\mathbb C}$ and the fiber derivatives of $F^2$ at $(o;v)$. After evaluating the base derivatives at $o$, we use the isomorphism $d\pi_e:\mathfrak m^{1,0}\to T_o^{1,0}M$ to regard the resulting expressions as functions of $v\in\mathfrak m^{1,0}\setminus\{0\}$.

Subtracting \eqref{eq:3.54} from \eqref{eq:3.53} and using the definition \eqref{eq:3.44} yields
\begin{align*}
\mathfrak R^i_{\;lk\bar j}&=\frac12c^\alpha_{\;\bar jk} c^i_{\;\alpha l}-\left(-\frac12c^\alpha_{\;\bar jk} c^i_{\;\alpha l}\right)=c^\alpha_{\;\bar jk}c^i_{\;\alpha l}.
\end{align*}

Finally, write $u=u^k\mathcal E_k$, $w=w^j\mathcal E_j$ and $y=y^l\mathcal E_l$. The definitions of the structure constants of $\mathfrak{g}^{\mathbb{C}}$ give
$$
[\sigma(w),u]_{\mathfrak h^{\mathbb C}}=\bar w^j u^k c^\alpha_{\;\bar jk}\mathcal F_\alpha,\qquad [\mathcal F_\alpha,\mathcal E_l]_{\mathfrak m^{1,0}}=c^i_{\;\alpha l}\mathcal E_i.
$$
Consequently,
$$
\mathfrak R_v(u,\sigma(w))y=\mathfrak R^i_{\;lk\bar j}\,y^l u^k\bar w^j\mathcal E_i=\bar w^j u^k y^l c^\alpha_{\;\bar jk}c^i_{\;\alpha l}\mathcal E_i=[[\sigma(w),u]_{\mathfrak h^{\mathbb C}},y]_{\mathfrak m^{1,0}}.
$$
This proves \eqref{eq:3.43} for arbitrary $u,w,y$.
\end{proof}

Taking $u=w$ and $y=v$ in Proposition~\ref{pro:3.17} recovers the contracted identity
\begin{equation}\label{eq:3.55}
\mathfrak R_v(w,\sigma(w))v
=[[\sigma(w),w]_{\mathfrak h^{\mathbb C}},v]_{\mathfrak m^{1,0}}.
\end{equation}

\begin{theorem}\label{th:3.18}
Let $(M=G/H,J)$ be a reductive complex homogeneous space with reductive decomposition $\mathfrak g=\mathfrak h\oplus\mathfrak m$ and $\operatorname{Ad}(H)\mathfrak m\subset\mathfrak m$. Let $F:T^{1,0}M\to[0,\infty)$ be a $G$-invariant complex Finsler metric on $M$. Then for the decomposition of the contracted horizontal $(1,1)$-curvature of the complex Rund connection $\mathbf{R}^i_{\;k\bar{j}}=\mathfrak{R}^i_{\;k\bar{j}}+\mathcal{R}^i_{\;k\bar{j}}$ as in \eqref{eq:3.12}, we have
\begin{align}
g_v(\mathbf{R}_v(w,\sigma(w))v, v) =g_v([[\sigma(w),w]_{\mathfrak{h}^{\mathbb{C}}},v]_{\mathfrak{m}^{1,0}},v)   +g_v(\mathcal{R}_v(w,\sigma(w))v,v),\label{eq:3.56}
\end{align}
where $v\in\mathfrak{m}^{1,0}\backslash\{0\}$ is the radial vector with variable direction, $w\in\mathfrak{m}^{1,0}\backslash\{0\}$ is an arbitrary vector and $\mathcal{R}_v(w,\sigma(w))v$ is given by \eqref{eq:3.41}.
\end{theorem}

\begin{proof}
Combining \eqref{eq:3.55} with \eqref{eq:3.12} and taking the $g_v$-inner product with $v$ yields \eqref{eq:3.56}.
\end{proof}

\begin{corollary}
Under the same assumptions as in Theorem \ref{th:3.18}, let $v\in\mathfrak{m}^{1,0}\backslash\{0\}$ be the radial vector with variable direction and let $w\in\mathfrak{m}^{1,0}\backslash\{0\}$ be an arbitrary vector. The holomorphic bisectional curvature in the directions $v$ and $w$ is given by
\begin{equation*}
B(v, w) = \frac{g_v\bigl(\mathbf{R}_v(w, \sigma(w))v, v\bigr)}{g_v(v, v)\, g_v(w, w)},
\end{equation*}
where $g_v(\mathbf{R}_v(w,\sigma(w))v,v)$ is given by \eqref{eq:3.56}.
\end{corollary}

\begin{corollary}
Under the assumptions of Theorem \ref{th:3.18}, let $v\in\mathfrak{m}^{1,0}\backslash\{0\}$ be the radial vector with variable direction. The holomorphic sectional curvature in the direction $v$ with  $g_v(v,v)=1$ is given by
\begin{equation*}
K(v) = 2g_v([[\sigma(v),v]_{\mathfrak{h}^{\mathbb{C}}},v]_{\mathfrak{m}^{1,0}},v)       +2g_v(\mathcal{R}_v(v,\sigma(v))v,v),
\end{equation*}
where $\mathcal{R}_v(v,\sigma(v))v$ is given by \eqref{eq:3.42}.
\end{corollary}

\section{Complex Finsler symmetric spaces}\label{sec:4}
\hspace{0.4cm}
We apply the algebraic criteria and complex Rund curvature formulas of Section~\ref{sec:3} to globally symmetric complex Finsler spaces. Their complex symmetric coset descriptions yield vanishing relations for the structure constants, from which the K\"ahler--Berwald property and the holomorphic curvature formulas follow.

\subsection{Complex symmetric coset spaces and Lie algebra relations}
\begin{definition}[cf. {\cite{DengHou2006}}]
Let $(M,J)$ be a connected complex manifold equipped with a strongly convex complex Finsler metric $F$. We call $(M,J,F)$ a \emph{globally symmetric complex Finsler space} if, for every $z\in M$, there exists a globally defined holomorphic isometry $s_z:M\to M$ satisfying
$$
s_z(z)=z,\qquad s_z^2=\mathrm{Id}_M,\qquad (ds_z)_z=-\mathrm{Id}_{T_zM}.
$$
The map $s_z$ is called a \emph{symmetry at $z$}. Its differential $ds_z$ acts on complex tangent vectors by complex linear extension.
\end{definition}

\begin{remark}
For a holomorphic involution fixing $z$, the condition $(ds_z)_z=-\mathrm{Id}_{T_zM}$ is equivalent to $z$ being an isolated fixed point. We  need only strongly pseudoconvexity of complex Finsler metrics in this paper.
\end{remark}

Now, let's recall the local properties of the geodesic equation of the complex Rund connection, which will be used to prove the transitivity in Theorem \ref{th:4.2}.

Let $(M,J,F)$ be a connected strongly pseudoconvex complex Finsler manifold, and let $\overset{\textbf{CR}}{\nabla}$ denote its complex Rund connection. Choose  holomorphic coordinates $(z^i)$ on $U\subset M$. Writing $v=v^i\partial/\partial z^i$ gives the induced holomorphic coordinates $(z^i;v^i)$ on $\tilde\pi^{-1}(U)\subset\tilde M$. Let $\overset{\mathcal{O}}{\Gamma}{}^i_{\;jk}(z;v)$ denote the coefficients of $\overset{\textbf{CR}}{\nabla}$ relative to the pulled-back holomorphic coordinate frame $\{\tilde\pi^*(\partial/\partial z^i)\}$. By Theorem~3.1 and Corollary~3.1 in \cite{LuoXiaoZhong2025}, in the pulled-back coordinate frame,
the corresponding connection 1-forms are $
\overset{\mathcal{O}}{\Gamma}{}^i_{\;jk}(z;v)\,dz^k,
$
here $dz^k$ also denotes the pullback of the corresponding coordinate $1$-form to $\tilde M$.

For a regular real curve $\gamma:[0,1]\to U$, write $z^i(t)=z^i(\gamma(t))$. The $(1,0)$ part of its velocity is
$$
\eta(t)=\frac12\bigl(\dot\gamma(t)-\sqrt{-1}J\dot\gamma(t)\bigr) =\dot z^i(t)\frac{\partial}{\partial z^i}\Big|_{\gamma(t)}.
$$
The corresponding geodesic equation is
\begin{equation}\label{eq:4.1}
\ddot z^i+\overset{\mathcal{O}}{\Gamma}{}^i_{\;jk}(z;\dot z)\dot z^j\dot z^k=0.
\end{equation}
By the complex homogeneity of $F$ and the cited coefficient formulas, \eqref{eq:4.1} is reversible.

For a real initial velocity $X\in T_{z_0}M$, let $\gamma_{z_0,X}$ be the solution with initial position $z_0$ and velocity $X$, and define
$$
\operatorname{Exp}^{\textbf{CR}}_{z_0}(X):=\gamma_{z_0,X}(1)
$$
whenever it is defined, taking $\gamma_{z_0,0}\equiv z_0$. The following lemma states its local properties.

\begin{lemma}\label{lm:4.1}
For every $z_0\in M$, there exists a neighborhood $U_{z_0}$ of $0$ in $T_{z_0}M$ such that 
$$\operatorname{Exp}^{\mathbf{CR}}_{z_0}:U_{z_0}\to \operatorname{Exp}^{\mathbf{CR}}_{z_0}(U_{z_0})$$ is a $C^1$ diffeomorphism onto a neighborhood of $z_0$ satisfying $(d\operatorname{Exp}^{\mathbf{CR}}_{z_0})_0=\mathrm{Id}_{T_{z_0}M}.$

Every global holomorphic isometry $f$ satisfies
$$
f(\operatorname{Exp}^{\mathbf{CR}}_{z_0}(X))=\operatorname{Exp}^{\mathbf{CR}}_{f(z_0)}(df_{z_0}(X))
$$
for all sufficiently small $X\in T_{z_0}M$, and is uniquely determined by its value and differential at one point.
\end{lemma}

\begin{proof}
In local holomorphic coordinates, the velocity-dependent term in \eqref{eq:4.1} is
\begin{equation}\label{eq:4.2a}
\overset{\mathcal{O}}{\Gamma}{}^i_{\;jk}(z;\dot z)\,\dot z^j\dot z^k.
\end{equation}
Here $\dot z=(\dot z^1,\dots,\dot z^m)$ is the coordinate velocity of a regular curve $\gamma:[0,1]\to U$, and $|\dot z|$ denotes its Euclidean norm. The displayed expression \eqref{eq:4.2a} is smooth for $\dot z\ne0$ and positively homogeneous of degree two in $\dot z$. For each $i$, the quantity in \eqref{eq:4.2a} and its first derivatives with respect to $z^l$ and $\bar z^l$, with $\dot z$ held fixed, are $O(|\dot z|^2)$ as $\dot z\to0$. Its first derivatives with respect to $\dot z^l$ and $\dot{\bar z}^l$ are $O(|\dot z|)$. These estimates hold locally uniformly in the base point. Thus, extending the quantity in \eqref{eq:4.2a} by zero at $\dot z=0$ gives a $C^1$ function of $(z,\dot z)$. The $C^1$ dependence of solutions on initial data implies that $\operatorname{Exp}^{\textbf{CR}}_{z_0}$ is $C^1$ near $0$. Linearizing the geodesic equation along the constant solution gives
$$
(d\operatorname{Exp}^{\textbf{CR}}_{z_0})_0=\mathrm{Id}_{T_{z_0}M}.
$$
The local diffeomorphism assertion follows from the inverse function theorem; see also \cite[Section~5.3, pp.~126--127]{bao2012introduction} for the corresponding exponential map argument.

The remaining assertions follow from the invariance of $\overset{\textbf{CR}}{\nabla}$ under holomorphic isometries, uniqueness of geodesics, and the identity theorem.
\end{proof}

The following theorem gives a homogeneous description of globally symmetric complex Finsler spaces in terms of complex symmetric coset spaces, whose definition is given below.
\begin{theorem}\label{th:4.2}
Let $(M,J,F)$ be a globally symmetric complex Finsler space. Then there exists a complex symmetric coset space $(G/H,\widetilde J,\Theta)$ equipped with a $G$-invariant complex Finsler metric $\widetilde F$ such that $(M,J,F)$ is holomorphically isometric to $(G/H,\widetilde J,\widetilde F)$. Conversely, every complex symmetric coset space $(G/H,J,\Theta)$ equipped with a $G$-invariant complex Finsler metric $F$ is a globally symmetric complex Finsler space.
\end{theorem}
To prove Theorem~\ref{th:4.2}, we first recall the symmetric coset space construction and introduce its complex counterpart.

Let $G$ be a connected real Lie group and $H$ a closed subgroup. The coset space $G/H$ is called symmetric if $G$ admits an involutive automorphism $\Theta$ such that
$$
(G_\Theta)_o\subset H\subset G_\Theta, \qquad G_\Theta=\{g\in G:\Theta(g)=g\},
$$
where $(G_\Theta)_o$ denotes the identity component of $G_\Theta$. The symmetry at $o=eH$ is
$$
\rho_o(gH)=\Theta(g)H,
$$
and the symmetry at $z=gH$ is $\rho_z=g\circ\rho_o\circ g^{-1}$, where elements of $G$ act by left translation.

Put $\theta:=d\Theta_e$. The $+1$ and $-1$ eigenspaces of $\theta$ in $\mathfrak{g}$ are $\mathfrak h$ and $\mathfrak m$, respectively, giving the canonical decomposition $\mathfrak g=\mathfrak h\oplus\mathfrak m$.
\begin{proposition}[{\cite[Chapter XI, Proposition 2.1]{FoundationOfDG2}}]\label{pro:4.2}
The canonical decomposition of a symmetric Lie algebra satisfies
$$
[\mathfrak h,\mathfrak h]\subset\mathfrak h,\qquad [\mathfrak h,\mathfrak m]\subset\mathfrak m,\qquad [\mathfrak m,\mathfrak m]\subset\mathfrak h.
$$
\end{proposition}
For every $h\in H\subset G_\Theta$, we have $\theta\circ\operatorname{Ad}(h)=\operatorname{Ad}(h)\circ\theta$. Thus $\operatorname{Ad}(H)$ preserves $\mathfrak m$ and a symmetric coset space is a reductive homogeneous space.

\begin{definition}
A complex symmetric coset space is a triple $(G/H,J,\Theta)$ consisting of an even-dimensional symmetric coset space $(G/H,\Theta)$ and a $G$-invariant almost complex structure $J$ on $G/H$.
\end{definition}
The $G$-invariant almost complex structure $J$ is automatically integrable: by Proposition~\ref{pro:4.2}, all projected brackets in the integrability condition of Section~\ref{sec:2.1} vanish. Moreover, under the identification $d\pi:\mathfrak m\to T_o(G/H)$,
\begin{equation*}
(d\rho_o)_o\circ d\pi=d\pi\circ\theta,\qquad (d\rho_o)_o=-\mathrm{Id}_{T_o(G/H)}.
\end{equation*}
Hence $(d\rho_o)_o$ commutes with $J_o$. The identity $\rho_o\circ g=\Theta(g)\circ\rho_o$ and the $G$-invariance of $J$ then show that $\rho_o$ is holomorphic on $G/H$, as are all its conjugates $\rho_z$.

\begin{proof}[Proof of Theorem~\ref{th:4.2}]
We follow the arguments of Deng--Hou \cite[Theorems~2.1 and~2.4]{DengHou2007}, with the necessary adaptations to the complex Finsler setting.

Let $\mathcal A=\mathcal A(M,J,F)$ be the group of holomorphic isometries, which is a Lie transformation group \cite[Theorem~1.5]{Spiro2001}. Fix $z_0\in M$ and choose the neighborhood $U_{z_0}\subset T_{z_0}M$ from Lemma~\ref{lm:4.1}. For $z_1=\operatorname{Exp}^{\textbf{CR}}_{z_0}(X)$ with $X\in U_{z_0}$, let $\gamma(t)=\gamma_{z_0,X}(t)$, $0\le t\le1$, and put $z_*=\gamma(1/2)$. By Lemma~\ref{lm:4.1} and the reversibility of \eqref{eq:4.1}, both $s_{z_*}\circ\gamma$ and $t\mapsto\gamma(1-t)$ are geodesics. At $t=1/2$, they have the same position $z_*$ and the same velocity $-\dot\gamma(1/2)$. Uniqueness for \eqref{eq:4.1} therefore gives
$$
s_{z_*}(\gamma(t))=\gamma(1-t).
$$
In particular, $s_{z_*}(z_0)=z_1$. By Lemma~\ref{lm:4.1}, points sufficiently close to any $z\in M$ can be joined to $z$ by geodesic segments. Since $M$ is connected, any two points can therefore be joined by a finite chain of geodesic segments. Applying the preceding midpoint argument to each segment, we obtain a finite composition of holomorphic symmetries sending one point to the other. Hence $\mathcal A$ acts transitively on $M$.

Since $M$ is connected, the identity component $G=\mathcal A^0$ also acts transitively. Fix $o\in M$ and let
$$
H=\{g\in G:g(o)=o\}.
$$
Then $H$ is closed, and the orbit map
$$
\Psi:G/H\longrightarrow M,\qquad \Psi(gH)=g(o),
$$
is a diffeomorphism. Set $\widetilde J:=\Psi^*J$ and define
$$
\widetilde F(z;v) :=F\bigl(\Psi(z);d\Psi_z(v)\bigr).
$$
Then $\widetilde J$ and $\widetilde F$ are $G$-invariant, and $\Psi$ is a holomorphic isometry.

Conjugation by $s_o$ preserves $G=\mathcal A^0$ and defines an involutive automorphism
$$
\Theta(g):=s_o\circ g\circ s_o^{-1},\qquad g\in G.
$$
For $h\in H$, the map $h\circ s_o\circ h^{-1}$ fixes $o$ and has differential $-\mathrm{Id}_{T_oM}$ there. By Lemma~\ref{lm:4.1}, it equals $s_o$. Thus $h$ commutes with $s_o$, so $\Theta(h)=h$ and $H\subset G_\Theta$.

Conversely, for $g\in(G_\Theta)_o$, commutation with $s_o$ gives $s_o(g(o))=g(s_o(o))=g(o)$. Hence the orbit $(G_\Theta)_o\cdot o$ is a connected subset of the fixed point set of $s_o$ containing $o$. Since $o$ is an isolated fixed point, this orbit is $\{o\}$, and therefore $(G_\Theta)_o\subset H$. We have proved
$$
(G_\Theta)_o\subset H\subset G_\Theta.
$$
It follows that $(G/H,\Theta)$ is a symmetric coset space. Together with the pulled-back $G$-invariant complex structure $\widetilde J$, this gives the required complex symmetric coset space $(G/H,\widetilde J,\Theta)$. Identifying $M$ with $G/H$ via $\Psi$, we henceforth write $J$ and $F$ in place of $\widetilde J$ and $\widetilde F$.

Conversely, let $(G/H,J,\Theta)$ be a complex symmetric coset space equipped with a
$G$-invariant complex Finsler metric $F$. As shown above, the induced symmetry
$$
\rho_o(gH)=\Theta(g)H
$$
is a holomorphic involution satisfying $(d\rho_o)_o=-\mathrm{Id}$. To verify that it preserves $F$, take $z=gH$ and write $w=dg_o(v)$ with $v\in T_o^{1,0}(G/H)$. Using $\rho_o\circ g=\Theta(g)\circ\rho_o$ and the $G$-invariance of $F$, we obtain
\begin{align*}
F\bigl(\rho_o(z);(d\rho_o)_z w\bigr)
=F\bigl(\Theta(g)H;d(\Theta(g))_o(-v)\bigr)
%&=F(o;-v)\\
%&=F(o;v)\\
=F(z;w),
\end{align*}
where the third equality follows from complex homogeneity. Thus $\rho_o$ is a global holomorphic isometry.

For each $z=gH$, its conjugate
$
\rho_z=g\circ\rho_o\circ g^{-1}
$
is also a global holomorphic isometry and satisfies
$$
\rho_z(z)=z,\qquad \rho_z^2=\mathrm{Id},\qquad (d\rho_z)_z=-\mathrm{Id}.
$$
Hence $(G/H,J,F)$ is a globally symmetric complex Finsler space.
\end{proof}

We now identify $M$ with $G/H$ and use the canonical reductive decomposition above. As in Section~\ref{sec:2.1}, let $I:\mathfrak m\to\mathfrak m$ correspond to $J_o$ and write
$\mathfrak g^{\mathbb C}=\mathfrak h^{\mathbb C}\oplus\mathfrak m^{\mathbb C}$ and $\mathfrak m^{\mathbb C}=\mathfrak m^{1,0}\oplus\mathfrak m^{0,1}.$
Complexifying $[\mathfrak m,\mathfrak m]\subset\mathfrak h$ immediately gives
\begin{equation}\label{eq:4.2}
[\mathfrak m^{1,0},\mathfrak m^{1,0}]\subset\mathfrak h^{\mathbb C},\qquad [\mathfrak m^{1,0},\mathfrak m^{0,1}]\subset\mathfrak h^{\mathbb C}.
\end{equation}
Moreover, the $\operatorname{Ad}(H)$-invariance of $I$ implies $[\mathfrak h^{\mathbb C},\mathfrak m^{1,0}]\subset\mathfrak m^{1,0}$ and its conjugate relation. For bases $\{\mathcal E_i\}$ of $\mathfrak m^{1,0}$ and $\{\mathcal F_\alpha\}$ of $\mathfrak h^{\mathbb C}$, \eqref{eq:4.2} says that all $\mathfrak m^{\mathbb C}$-components of $[\mathcal E_j,\mathcal E_k]$ and $[\mathcal E_j,\sigma(\mathcal E_k)]$ vanish. In the notation of Section~\ref{sec:2.1}, this is precisely
\begin{equation}\label{eq:4.3}
c^i_{\;jk}=0,\qquad c^i_{\;j\bar k}=0,\qquad c^{\bar i}_{\;j\bar k}=0.
\end{equation}

\subsection{Rigidity theorem}
\hspace{0.4cm}
The vanishing relations \eqref{eq:4.3} simplify the K\"ahler-Finsler criterion and the connection coefficients of Section~\ref{sec:3}, yielding the following rigidity theorem.

\begin{theorem}\label{th:4.5}
Let $(G/H,J,\Theta)$ be a complex symmetric coset space equipped with a $G$-invariant strongly pseudoconvex complex Finsler metric $F$. Then $F$ is a K\"ahler--Berwald metric.
\end{theorem}

\begin{proof}
Let $\mathfrak g=\mathfrak h\oplus\mathfrak m$ be the canonical decomposition associated with $\theta=d\Theta_e$. By Theorem~\ref{th:3.16}, $F$ is K\"ahler--Finsler if and only if
\begin{align}
c_{k\bar q o}-c_{\bar q ok}-c_{o\bar q k}-g_{\bar q\bar j}\,c^{\bar j}_{\ \bar o k}-C_{k\bar q l}\,\mathbf{N}^l_{\;o}+C_{k\bar q\bar j}\,c^{\bar j}_{\ \bar o o}=0.\label{eq:4.4}
\end{align}
Upon substituting \eqref{eq:4.3} into \eqref{eq:4.4}, all terms involving the structure constants $c^i_{\;jk}$, $c^i_{\;j\bar k}$, $c^{\bar i}_{\;j\bar k}$ vanish. Moreover, the connection coefficients $\mathbf{N}^j_{\;k}$ are defined via $\mathbf{N}^j_{\;k}=c^{\;\,j}_{o\;k}+g^{\bar q j}g_{\bar q\bar l}\,c^{\bar l}_{\;\bar o k}$, which also vanish by \eqref{eq:4.3}. Hence every term in \eqref{eq:4.4} is zero, and $F$ is a K\"ahler--Finsler metric.

The connection coefficients of the complex Rund connection associated with $F$ are given by
\begin{equation}\label{eq:4.5}
\mathbf{\Gamma}^i_{\;jk}(o;v)=c^{\;i}_{j\;k}-C^{\;i}_{j\;l}\,\mathbf{N}^l_{\;k}+C^{\;i}_{j\;\bar l}\,c^{\bar l}_{\;\bar o k}.
\end{equation}
Substituting \eqref{eq:4.3} into \eqref{eq:4.5} yields $\mathbf{\Gamma}^i_{\;jk}(o;v)=0$ for every nonzero $v\in\mathfrak m^{1,0}$. Consequently,
$$
\frac{\partial}{\partial v^l}\mathbf{\Gamma}^i_{\;jk}(o;v)=\frac{\partial}{\partial\bar v^l}\mathbf{\Gamma}^i_{\;jk}(o;v)=0.
$$
Since the complex Rund connection is $G$-invariant, its fiber independence at the origin propagates to every point of $G/H$. Hence $F$ is a complex Berwald metric. Since $F$ is also K\"ahler--Finsler, it is a K\"ahler--Berwald metric.
\end{proof}

This together with Theorem~\ref{th:4.2} show that every globally symmetric complex Finsler space is a K\"ahler--Berwald space.

\begin{remark}
The existence of a $G$-invariant K\"ahler--Berwald metric on a reductive complex homogeneous space $(G/H,J)$ does not imply  $G/H$ being a  complex symmetric coset space.
\end{remark}

\begin{example}
Let
$$
G=SU(3)\times\mathbb R^2,\qquad H=T^2\times\{(0,0)\},
$$
where $T^2$ is the diagonal maximal torus of $SU(3)$. With respect to the standard $G$-invariant complex structure, choose the reductive complement $\mathfrak m$ such that
$$
\mathfrak m^{1,0}=\operatorname{span}_{\mathbb C}\{\mathcal E_1,\mathcal E_2,\mathcal E_3,\mathcal E_4\},\qquad \mathcal E_1=E_{12},\quad \mathcal E_2=E_{23},\quad \mathcal E_3=E_{13},\quad \mathcal E_4=e_1-\sqrt{-1}e_2.
$$
Here, $e_1,e_2$ form the standard basis of  $\mathbb R^2$, and $E_{ij}$ denotes the matrix unit in $\mathfrak{sl}_3(\mathbb{C})$.  On $\mathfrak{su}(3)^{\mathbb C}\cong\mathfrak{sl}_3(\mathbb C)$, the conjugation defining the real form $\mathfrak{su}(3)$ is $\sigma(X)=-\overline{X}^{\mathsf T}$ for $X\in \mathfrak{su}(3)^{\mathbb C}$, so that $\sigma(E_{ij})=-E_{ji}$. In particular, $[\mathcal E_1,\mathcal E_2]=\mathcal E_3$.

For $v=v^i\mathcal E_i$ and each integer $k\geq1$, define
$$
F_k^2(v)=\mathsf{v}+\mathsf{w}+\sqrt[k]{\mathsf{v}^{2k}+\mathsf{w}^{2k}},
$$
where $\mathsf{v} =|v^1|^2+|v^2|^2+2|v^3|^2$ and $\mathsf{w} =|v^4|^2$ and the root is taken to be nonnegative. The isotropy action of $H$ on $\mathfrak{m}^{1,0}$ multiplies $\mathcal E_1,\mathcal E_2,\mathcal E_3$ by unit complex numbers and fixes $\mathcal E_4$. Hence $F_k$ is $\operatorname{Ad}_{\mathfrak m^{1,0}}(H)$-invariant and induces a $G$-invariant complex Finsler metric on $G/H$.

The quadratic form $\mathsf{v}$ is, up to an overall normalization,  the squared norm at the origin of the standard $SU(3)$-invariant K\"ahler--Einstein metric on $SU(3)/T^2$, while $\mathsf{w}$ corresponds to the standard flat Hermitian metric on $\mathbb C$. By Zhong's product construction \cite[Theorem 10]{Zhong2023}, for every $k\geq 1$, the induced metric is strongly convex, non-Hermitian and K\"ahler--Berwald. Every $\operatorname{Ad}(H)$-invariant complement $\mathfrak{m}$ of $\mathfrak h$ has  complexification $\mathfrak{m}^{\mathbb{C}}$ containing the root spaces $\mathbb{C}E_{12},\mathbb{C}E_{23},\mathbb{C}E_{13}$. Indeed, these are nonzero $\mbox{Ad}(H)$-weight spaces, so they cannot lie in $\mathfrak{h}^{\mathbb{C}}$. If $\mathfrak{m}$ defined a symmetric decomposition, one would require $[\mathfrak{m},\mathfrak{m}]\subset\mathfrak{h}$, and hence, after complexification, $[\mathfrak{m}^{\mathbb{C}},\mathfrak{m}^{\mathbb{C}}]\subset\mathfrak{h}^{\mathbb{C}}$. But $E_{12},E_{23}\in\mathfrak{m}^{\mathbb{C}}$ and
$$
[E_{12},E_{23}]=E_{13}\notin\mathfrak h^{\mathbb C},
$$
which is a contradiction. Therefore no such $\mbox{Ad}(H)$-invariant complement can define a symmetric decomposition. Thus the homogeneous pair $(G,H)$ admits no symmetric decomposition, and with respect to the standard homogeneous complex structure,  $G/H$ is not a complex symmetric coset space.
\end{example}

In the symmetric case, the holomorphic bisectional curvature and holomorphic sectional curvature formulas can be greatly simplified as follows.
\begin{theorem}\label{th:4.7}
Let $(G/H,J,\Theta)$ be a complex symmetric coset space equipped with a $G$-invariant strongly pseudoconvex complex Finsler metric $F$, and let $\mathfrak g=\mathfrak h\oplus\mathfrak m$ be the canonical decomposition. Let $v\in\mathfrak{m}^{1,0}\backslash\{0\}$ and $w\in\mathfrak{m}^{1,0}\backslash\{0\}$ such that $g_v(v,v)=g_v(w,w)=1$. Then the holomorphic bisectional curvature in the directions  $v$ and $w$,  and the holomorphic sectional curvature in the direction $v$ are respectively given  by
\begin{align*}
B(v,w)=g_v\bigl([[\sigma(w),w]_{\mathfrak h^{\mathbb C}},v]_{\mathfrak m^{1,0}},v\bigr)\quad\mbox{and}\quad
K(v)=2g_v\bigl([[\sigma(v),v]_{\mathfrak h^{\mathbb C}},v]_{\mathfrak m^{1,0}},v\bigr).
\end{align*}
\end{theorem}

\begin{proof}
Since $[\mathfrak m,\mathfrak m]\subset\mathfrak h$, all $\mathfrak m^{\mathbb C}$-components of brackets between elements of $\mathfrak m^{\mathbb C}$ vanish. Equation~\eqref{eq:3.35} therefore gives $\mathcal N_v=0$ for every nonzero $v$. Substituting these vanishing relations into \eqref{eq:3.41} yields
$$
\mathcal R_v(w,\sigma(w))v=0.
$$
The formula for $B(v,w)$ now follows from Theorem~\ref{th:3.18} and the definition of holomorphic bisectional curvature. Taking $w=v$ and using $K(v)=2B(v,v)$ gives the formula for $K(v)$.
\end{proof}

\vskip0.4cm
\noindent
\textbf{Acknowledgements.}\ {\small C. Zhong was supported by the National Natural Science Foundation of China (No. 12471080) and the Key Program of the Natural Science Foundation of Fujian Province (No. 2026J002007)}.

\end{document}